\documentclass[11pt,a4paper]{article}

\usepackage[english]{babel}
\usepackage{amsmath}
\usepackage{amssymb,amsfonts,a4}
\usepackage{eucal,bbm,stmaryrd}
\usepackage{psfrag,graphicx}

\begin{document}

\numberwithin{equation}{section}

\newtheorem{satz}{Theorem}
\newtheorem{defi}{Definition}
\newtheorem{lem}{Lemma}

%Vorbereitungen

\newcommand{\Bew}{\noindent{\bf Proof: }}
\newcommand{\qed}{\hspace*{3em} \hfill{$\square$}}
\newcommand{\weg}{\setminus}
\newcommand{\gdw}{\Leftrightarrow}
\newcommand{\sgdw}{\Leftrightarrow}
\newcommand{\nach}{\longrightarrow}
\newcommand{\alle}{\,\forall\,}
\newcommand{\gibt}{\,\exists\,}
\newcommand{\Bed}{\,|\,}

\newcommand{\Z}{\mathbbm{Z}}
\newcommand{\R}{\mathbbm{R}}
\newcommand{\baR}{\overline{\R}}
\newcommand{\N}{\mathbbm{N}}

\newcommand{\p}{\mathcal{P}}
\newcommand{\prob}{\mathcal{P}_1}

\newcommand{\opB}{B_2}

%Raeume und Gibbsmasse

\newcommand{\X}{\mathcal{X}}
\newcommand{\XX}{\mathcal{X}_0}
\newcommand{\bX}{\bar{X}}
\newcommand{\La}{\Lambda}
\newcommand{\Lan}{\La_n}

\newcommand{\Bo}{\mathcal{B}}
\newcommand{\B}{\mathcal{B}}
\newcommand{\F}{\mathcal{F}}

\newcommand{\la}{\lambda}
\newcommand{\si}{\sigma}

\newcommand{\einh}{e_1}
\newcommand{\lear}{\le_{\einh}}
\newcommand{\Ga}{\Gamma}
\newcommand{\Gar}{\Ga_r}
\newcommand{\Gas}{\Ga_s}

\newcommand{\bs}{\bar{s}}

\newcommand{\E}{\mathcal{E}}

\newcommand{\gmu}{\gamma}
\newcommand{\G}{\mathcal{G}}

%Hauptteil

\newcommand{\gglatt}{\varphi}

\newcommand{\en}{E_{n}}

\newcommand{\ah}{\vec{a}} 
\newcommand{\tauh}{\vec{\tau}} 

\newcommand{\Nk}{\N_k}

\newcommand{\A}{\mathcal{A}}
\newcommand{\K}{\mathcal{K}}
\newcommand{\Cinf}{\mathcal{C}^{\infty}}

\newcommand{\Gn}{G_n}
\newcommand{\Gna}{G'_n}
\newcommand{\Gnb}{G''_n}

\newcommand{\tn}{\tau_n}
\newcommand{\Taun}{\tau^q_n}
\newcommand{\bn}{\bar{n}}
\newcommand{\bR}{\bar{R}}
\newcommand{\bU}{\bar{U}}
\newcommand{\tU}{\tilde{U}}
\newcommand{\tiu}{\tilde{u}}

\newcommand{\pin}{\pi_n}
\newcommand{\pinh}{\pi'_n}

\newcommand{\mxt}{m_{x',t}}
\newcommand{\hxt}{h_{x',t}}

\newcommand{\bpsi}{\bar{\psi}}

\newcommand{\de}{\delta}
\newcommand{\ep}{\epsilon}
\newcommand{\ex}{c_\xi}
\newcommand{\bB}{\bar{B}}

\newcommand{\ph}{\varphi_n}
\newcommand{\bph}{\bar{\varphi}_n}
\newcommand{\leb}{\lambda^2}
\newcommand{\cK}{c_K}
\newcommand{\cpsi}{c_{\psi}}
\newcommand{\cu}{c_u}
\newcommand{\KU}{K^{U}}
\newcommand{\Ko}{K}
\newcommand{\Kop}{K'}
\newcommand{\Kopp}{K''}
\newcommand{\Kep}{K_{\ep}}
\newcommand{\Kepep}{K_{\ep,\ep}}
\newcommand{\fc}{f_{\Ko}}
\newcommand{\cf}{c_f}
\newcommand{\ft}{\tilde{f}}

\newcommand{\per}{\eta}
\newcommand{\akb}{A_{k,B,\per}}
\newcommand{\takb}{\tilde{A}_{k,B,\per}}

%%%

\newcommand{\anxb}{a_{n,X,B_+}}
\newcommand{\rnxb}{r_{n,X,B_+}}

\newcommand{\T}{\mathfrak{T}}
\newcommand{\Tn}{\T_n}
\newcommand{\Tnx}{\T_{n,X}}
\newcommand{\Tnb}{\T_{n,B}}

\newcommand{\tnxb}{t_{n,X,B}}
\newcommand{\Tnxb}{T_{n,X,B}}
\newcommand{\taunxb}{\tau_{n,X,B}}
\newcommand{\pnxb}{P_{n,X,B}}
\newcommand{\cnxb}{C_{n,X,B}}

\newcommand{\etn}{t_{n}}
\newcommand{\eTn}{T_{n}}
\newcommand{\etaun}{\tau_{n}}
\newcommand{\ela}{\Lambda}
\newcommand{\epn}{P_n}
\newcommand{\ecn}{C_n}

\newcommand{\tX}{\tilde{X}}
\newcommand{\tB}{\tilde{B}}
\newcommand{\tx}{\tilde{x}}

\newcommand{\tTn}{\tilde{\T}_n}
\newcommand{\tTnx}{\tilde{\T}_{n,\tX}}
\newcommand{\tTnb}{\tilde{\T}_{n,\tB}}

\newcommand{\ttnxb}{\tilde{t}_{n,\tX,\tB}}
\newcommand{\tTnxb}{\tilde{T}_{n,\tX,\tB}}
\newcommand{\ttaunxb}{\tilde{\tau}_{n,\tX,\tB}}
\newcommand{\tpnxb}{\tilde{P}_{n,\tX,\tB}}
\newcommand{\tcnxb}{\tilde{C}_{n,\tX,\tB}}

\newcommand{\ettn}{\tilde{t}_{n}}
\newcommand{\etTn}{\tilde{T}_{n}}
\newcommand{\ettaun}{\tilde{\tau}_{n}}
\newcommand{\etpn}{\tilde{P}_n}
\newcommand{\etcn}{\tilde{C}_n}

\newcommand{\iTn}{\bar{\T}_n}
\newcommand{\iTnx}{\bar{\T}_{n,X}}
\newcommand{\iTnB}{\bar{\T}_{n,B}}
\newcommand{\iTnxb}{\bar{T}_{n,X,B}}

%Hardcore

\newcommand{\anx}{a_{n,X}}
\newcommand{\rnx}{r_{n,X}}
\newcommand{\Kepn}{\Kep^n}
\newcommand{\hcny}{C_{X,\Kepn}}

\newcommand{\tnx}{t_{n,X}}
\newcommand{\taunx}{\tau_{n,X}}
\newcommand{\pnx}{P_{n,X}}
\newcommand{\Tnxn}{T_{n,X}}

\newcommand{\ttnx}{\tilde{t}_{n,\tX}}
\newcommand{\tTnxn}{\tilde{T}_{n,\tX}}
\newcommand{\ttaunx}{\tilde{\tau}_{n,\tX}}
\newcommand{\tpnx}{\tilde{P}_{n,\tX}}
\newcommand{\ak}{A_{k,\per}}
\newcommand{\tak}{\tilde{A}_{k,\per}}

%========================================================================

\begin{center}
{\bf \LARGE Translation-invariance of two-dimensional
 Gibbsian point processes}\\

\vspace{ 1 cm}
Thomas Richthammer\\
Mathematisches Institut der Universit\"at M\"unchen\\
Theresienstra\ss e 39, D-80333 M\"unchen\\
Email: Thomas.Richthammer@mathematik.uni-muenchen.de\\
Tel: +49 89 2180 4633\\
Fax: +49 89 2180 4043 
\end{center}

\bigskip

\begin{abstract}
The conservation of translation as a symmetry in two-dimensional systems 
with interaction is a classical subject of statistical mechanics. 
Here we establish such a result for Gibbsian particle systems 
with two-body interaction, where the interesting cases of singular, hard-core 
and discontinuous interaction are included. We start with the special 
case of pure hard core repulsion  in order to show 
how to treat hard cores in general.\\
 
Key words: Gibbsian point processes, Mermin-Wagner theorem, translation, 
hard-core potential, singular potential, pure hard core repulsion, 
percolation, superstability. 
\end{abstract}

%=========================================================================

\section{Introduction}

\begin{sloppypar}

Gibbsian processes were introduced by  R.~L.~Dobrushin (see \cite{D68}
and \cite{D70}), O.~E.~Lanford and D.~Ruelle (see \cite{LR}) as a
model for equilibrium states in statistical physics. (For general
results on Gibbs measures on a d-dimensional lattice we refer to the
books of H.-O.~Georgii \cite{G}, B.~Simon \cite{Sim} and Y.~G.~Sinai
\cite{Sin}, which cover a wide range of phenomena.) 
The first results concerned existence and uniqueness of
Gibbs measures and the structure of the set of Gibbs measures related
to a given potential. The question of uniqueness is of special
importance, as the non-uniqueness of Gibbs measures can be interpreted
as a certain type of phase transition occurring within the particle
system. A phase transition occurs whenever a symmetry of the
potential is broken, so it is natural to ask, under which
conditions symmetries are broken or conserved. The answer to this
question depends on the type of the symmetry (discrete or continuous),
the number of spatial dimensions and smoothness and decay conditions
on the potential (see \cite{G}, chapters 6.2, 8, 9 and 20). 
It turns out that the case of continuous
symmetries in two dimensions is especially interesting. The first
progress in this case was achieved by M.~D.~Mermin and H.~Wagner, who
showed for special two-dimensional lattice models that continuous internal
symmetries are conserved (\cite{MW} and \cite{M}).
In \cite{DS} R.~L.~Dobrushin and S.~B.~Shlosman established 
conservation of symmetries for more general 
potentials which satisfy smoothness and decay conditions, and
C.-E.~Pfister improved this in  \cite{P}. Later also continuum systems
were considered: 
S. Shlosman obtained results for  continuous internal symmetries (\cite{S}), 
while  J.~Fr\"ohlich and C.-E.~Pfister treated
the case of translation of point particles  (\cite{FP1} and \cite{FP2}).
All these results rely on the smoothness of the interaction, 
but in \cite{ISV} D.~Ioffe, S.~Shlosman and Y.~Velenik were able 
to relax this condition. Considering a lattice model they showed that 
continuous internal symmetries are conserved, whenever the interaction can be 
decomposed into a smooth part and a part which is small 
with respect to $L^1$-norm, using a perturbation expansion and  
percolation theory. 
We generalised this to a point particle setting (\cite{Ri1}).\\

Here we will investigate the conservation of translational symmetry for non-smooth, 
singular or hard-core potentials in a point particle setting. While we treat 
non-smoothness by generalising ideas used in \cite{Ri1}, we will give an
approach to singular potentials which is different from the one given in 
\cite{FP1} and \cite{FP2}. The advantage of our approach is that 
integrability condition (2.13) of \cite{FP2} is simplified and relaxed and 
the case of hard-core potentials can easily be included. Thus we 
are able to show the conservation of translational symmetry for the
pure hard core model, for example.\\ 

In Section \ref{secresult}  we will first confine ourselves to this special 
case of pure hard core repulsion. The corresponding result 
(Theorem \ref{purehardcore}) is of interest on its own and its proof shows 
how to deal with hard cores in the general case. For this general case we then 
define a suitable class of potentials (Definition \ref{defGapprox}), 
give some sufficient conditions for potentials to belong to that class 
(Lemmas \ref{lepotglatt} and~\ref{lepotchar}) and state the general result 
obtained (Theorem \ref{sym}). The precise setting is then given in 
Section \ref{secsetting}. The proofs of the lemmas from Sections \ref{secresult} 
and \ref{secsetting} are relegated to Section \ref{secleset}.
In Sections \ref{secproofcore} and \ref{secproofsym} we will 
give the proofs of Theorems \ref{purehardcore} and \ref{sym} respectively. 
The proofs of the corresponding lemmas are relegated to Sections
\ref{secleproofcore} and \ref{secleproofsym} respectively. 
In the proof of the general case arguments of the special case
have to be modified and refined by new concepts and ideas at several instances.  
So for sake of clarity we will repeat arguments from the proof of 
Theorem \ref{purehardcore} in the proof of Theorem \ref{sym} whenever 
necessary.\\ 

Acknowledgement: I would like to thank H.-O. Georgii 
for suggesting the problem and many helpful discussions and 
F.~Merkl for helpful comments.

%=====================================================================
%
%=====================================================================

\section{Result} \label{secresult}

We consider particles in the plane $\R^2$ without internal degrees of 
freedom. The chemical potential $-\log z$ of the system is given via an 
activity parameter $z>0$. The interaction between particles is modelled 
by a translation-invariant pair potential $U$, i.e. a measurable function 
\begin{displaymath}
U: \R^2 \; \to \; \baR \; := \; \R \cup \{\infty\},
\end{displaymath}
which is assumed to be symmetric in that $U(x)=U(-x)$ for all $x \in \R^2$. 
The potential of two particles $x_1,x_2 \in \R^2$ is then given by 
$U(x_1-x_2)$.\\

We first consider the particular case of pure hard core repulsion, where 
the size and the shape of the hard core are  given by a norm 
$|.|_h$ on $\R^2$. The corresponding pure hard-core potential $U_{hc}$ 
is defined by 
\[
U_{hc} (x) \, := \, 
\left \{
\begin{array}{ll}
\infty \quad &\text{for } \quad  |x|_h \le 1 \\
0 \quad &\text{for } \quad  |x|_h > 1. 
\end{array}
\right.
\] 
\begin{satz} \label{purehardcore}
Let $z>0$ be an activity  parameter, $|.|_h$ be a norm on $\R^2$ and 
$U_{hc}$ be the corresponding pure hard-core potential. 
Then every Gibbs measure corresponding to $U_{hc}$ and $z$ is translation-invariant.
\end{satz}   
The proof of Theorem \ref{purehardcore}, which is given in 
Section \ref{secproofcore}, will show how to deal with hard cores in the 
general case presented below.\\

In order to describe a class of potentials for which translational symmetry 
is conserved we will define important properties of sets, functions and potentials. 
A set $A \subset \R^2$ is called \emph{symmetric} if $A = -A$. 
We call $U$ a \emph{standard potential} if $U$ is a measurable, 
symmetric pair potential and its hard core 
\[
\KU \, := \, \{ U = +\infty\}
\] 
is bounded. Usually the hard core will be empty, $\{0\}$ or a disc, 
but in our setup we are able to treat fairly general hard cores.
For a given function $\psi:\R^2  \to \R_+$ we say that a standard
potential $U$ has  \emph{$\psi$-dominated} derivatives on the set $A$ 
if   
\[
\partial^2_i U(x+te_i) \, \le \, \psi(x)
\quad \text{ for all } x \in A, t \in [-1,1] \text{ s.t. }  
x+t e_i \in A
\]
for $i = 1,2$. Here $e_1 = (1,0)$, $e_2 = (0,1)$ and $\partial_i$ 
is the partial derivative in direction $e_i$. The above definition 
is meant to imply that these derivatives exist. In the context of 
$\psi$-domination we will use the notion of a \emph{decay function}, 
which is defined to satisfy
\[
\|\psi\| \, < \, \infty \quad \text{ and }  
\int \psi(x) |x|^2  dx  \, < \, \infty.
\]
This definition of course does not depend on the choice of norm $|.|$, 
but for sake of definiteness let $|.|$ be the maximum norm on $\R^2$.\\ 

If $U$ is a potential, $z$ is an activity parameter and $\XX$ is a
set of boundary conditions, we say that the triple $(U,z,\XX)$ is
admissible if all conditional Gibbs distributions corresponding to 
$U$ and $z$ with boundary condition taken from $\XX$ are well defined, 
see Definition \ref{defzul} in Section \ref{secgibbs}. 
Important examples are  the cases of superstable potentials with 
tempered boundary configurations and nonnegative potentials with arbitrary 
boundary conditions, see Section \ref{secadm}. For admissible 
$(U,z,\XX)$ the set of Gibbs measures $\G_{\XX}(U,z)$ corresponding to 
$U$ and $z$ with full weight on configurations in $\XX$ is a well defined 
object. Finally we need bounded correlations: For admissible $(U,z,\XX)$ 
we call $\xi \in \R$ a Ruelle bound if the correlation function of 
every Gibbs measure $\mu \in \G_{\XX}(U,z)$ is bounded by powers of $\xi$ 
in the sense of \eqref{ruellebd} in Section \ref{secgibbs}.
\begin{defi} \label{defGapprox}
Let $(U,z,\XX)$ be an admissible triple with Ruelle bound $\xi$, 
where $U:  \R^2\to \baR$ is a translation-invariant 
standard potential. We say that $U$ is smoothly approximable
if there is a decomposition of $U$ into a smooth part $\bU$ and
a small part $u$ in the following sense:
We have a symmetric, compact set $\Ko \supset \KU$,
a decay function  $\psi$ and measurable symmetric functions 
$\bU, u: \Ko^c \to \R$ such that 
\begin{equation} \label{Gapprox}
\begin{split}
&U =  \bU - u \text{ and }  u  \ge  0 \quad \text{ on } \Ko^c,\\
&\text{$\bU$ has  $\psi$-dominated derivatives on  $\Ko^c$},\\ 
&\int_{\Ko^c} \tiu (x)  |x|^2 \, dx \, < \, \infty  \, \text { and } \,
 \leb(\Ko \weg \KU) \, + \, \int_{\Ko^c} \tiu(x)\,dx \, < \, \frac 1 {z\xi},
\end{split}
\end{equation}
where $\tiu := 1 - e^{-u} \le  u \wedge 1$.
\end{defi}
The class of smoothly approximable standard potentials is a rich 
class of potentials. A smoothly approximable standard potential $U$ 
may have a singularity or a hard core at the origin, and the type of 
convergence into the singularity or the hard core is fairly arbitrary, 
as we have  not imposed any condition on $U$ in the set $\Ko \weg \KU$. 
For small activity $z$ the last condition of \eqref{Gapprox} holds for 
large sets $\Ko$, which relaxes the conditions on $U$. The small 
part $u$ of $U$ is not assumed to satisfy any regularity conditions, so that
$U$ doesn't have to be smooth or continuous.   
We note that Definition \ref{defGapprox} does not depend on the choice of 
the norm $|.|$. If we know a potential to be smooth outside of its hard core 
the above conditions simplify: 
\begin{lem} \label{lepotglatt}
Let $(U,z,\XX)$ be an admissible triple with Ruelle bound $\xi$, 
where $U: \R^2 \to \baR$ is a translation-invariant standard potential.
Suppose  we have a symmetric compact set $\Ko \supset \KU$ and a 
decay function  $\psi$ such that 
$U$ has $\psi$-dominated derivatives on $\Ko^c$ and 
$\leb(\Ko \weg \KU) < 1/(z \xi)$. 
Then $U$ is smoothly approximable.
\end{lem}
This is an immediate consequence of Definition \ref{defGapprox}.   
In the non-smooth case,  the following lemma
gives important examples of smoothly approximable potentials: 
\begin{lem} \label{lepotchar}
Let  $(U,z,\XX)$ be an admissible triple with Ruelle bound $\xi$, where 
$U: \R^2 \to \baR$ is a translation-invariant standard potential
such that $\KU$ is compact and $U$ is continuous in $(\KU)^c$.
Suppose we have a  decay function  $\psi$ and a compact set 
$\tilde{K} \subset \R^2$ such that one of the following properties holds:
\begin{enumerate}
\item[(a)] $U$ has $\psi$-dominated derivatives in $\tilde{K}^c$. 
\item[(b)] There is a standard potential $\tU \ge 0$ such that 
$|U| \le \tU$ in $\tilde{K}^c$, $\tU$  has $\psi$-dominated derivatives in 
$\tilde{K}^c$ and $\int_{\tilde{K}^c} \tU(x) |x|^2 dx < \infty$. 
\end{enumerate}
Then $U$ is smoothly approximable.
\end{lem}
For example, $(a)$ holds trivially when $U$ has finite range, 
and $(b)$ includes the case that there are  $\ep > 0$ and $k \ge 0$ 
such that $|U(x)| \le k / |x|^{4+\ep'}$ for large $|x|$. 
Our main result is now the following:  
\begin{satz} \label{sym}
Let $(U,z,\XX)$ be admissible with Ruelle bound, 
where $U: \R^2 \to \baR$ is a translation-invariant standard potential. 
If $U$ is smoothly approximable then 
every Gibbs measure $\mu \in \G_{\XX}(U,z)$ is translation-invariant. 
\end{satz}
For a generalisation of the above result to the case of particles
with inner degrees of freedom, i.e. Gibbsian systems of marked
particles, we refer to \cite{Ri2}. 

%======================================================================
%
%======================================================================

\section{Setting} \label{secsetting}

\subsection{State space}

We will use the notations $\N := \{0,1,\ldots\}$, 
$\R_+ := [0,\infty[$, $\bar{\R} := \R \cup \{+\infty\}$,  
\[
r_1 \vee r_2 \, := \, \max\{r_1,r_2\} \; \text{ and } \;     
r_1 \wedge r_2 \, := \, \min\{r_1,r_2\} \quad \text{ for } 
r_1,r_2 \in \R.
\]
On $\R^2$ we consider the maximum norm $|.|$ and the Euclidean norm $|.|_2$. 
For $\ep>0$ the \emph{$\ep$-enlargement} of a set $A \subset \R^2$ is 
defined by 
\[
A_{\ep}  \, := \,  \{x + x': x \in A, |x'|_2 < \ep\}. 
\]
The state space of a particle is the plane $\R^2$. The  Borel-$\si$-algebra 
$\B^2$ on $\R^2$ is induced by any norm on $\R^2$. Let $\B_b^2$ be the set 
of all bounded Borel sets and $\la^2$ be the Lebesgue measure on $(\R^2,\Bo^2)$. 
Integration with respect to this measure will be abbreviated by 
$dx := d\la^2(x)$. Often we consider the centred squares  
\[
\La_r \,:= \,[-r,r[^2 \, \subset \, \R^2 \qquad (r \in \R_+).
\]
We also want to consider bonds between particles. 
For a set $X$ we denote the set of all bonds in $X$ by 
\[
E(X)  \, :=  \, \{ A \subset X: \#A = 2\}.
\] 
A bond will be denoted by $xx' :=\{ x,x'\}$, where $x,x' \in X$ such that 
$x \neq x'$. For a bond set $B \subset E(X)$  $~(X,B)$ 
is an (undirected) graph, and we set 
\[
\begin{split}
x \stackrel{X,B}{\longleftrightarrow} x' \quad :\gdw \quad 
\gibt m \in \N, x_0,\ldots,x_m \in X: \, 
&x = x_0, x' = x_m,\\
&x_{i-1}x_{i} \in B \text{ for all } 1 \le i \le m. 
\end{split}
\]
This connectedness relation is an equivalence relation on $X$ 
whose equivalence classes are called the $B$-clusters of $X$. Let  
\[
C_{X,B}(x) \, :=\, \{ x' \in  X: x \stackrel{X,B}
                                {\longleftrightarrow} x' \}   
 \quad \text{ and } \quad 
 C_{X,B}(\La) \, := \,  \bigcup_{x' \in X \cap \La} C_{X,B}(x') 
\]
denote the $B$-clusters of a point $x$ and a set $\La$ respectively. 
Primarily we are interested in the case  $X = \R^2$. 
On the corresponding bond set  $E(\R^2)$ we consider the 
$\si$-algebra \[
\F_{E(\R^2)} \, := \, \{ \{ x_1x_2 \in E(\R^2):  
(x_1,x_2) \in M \} :  \, M \in (\B^2)^2 \}.
\]
Every symmetric function $u$ on $\R^2$ can be considered a function on 
$E(\R^2)$ via $u(xx') := u(x-x')$.

%=====================================================================

\subsection{Configuration space}

A set of particles $X \subset \R^2$ is called  
\[
\begin{array} {ll}
\text{finite } \quad &\text{ if } \#X < \infty, \qquad \text{ and }\\
\text{locally finite } \quad &\text{ if } \#(X \cap \La) < \infty
\text{ for all } \La \in \B^2_b,
\end{array}
\] 
where $\#$ denotes the cardinality of a set. The configuration space $\X$ 
of particles is defined as the set of all locally finite subsets of 
$\R^2$, and its elements are called configurations of particles. 
For $X,\bX \in \X$ let  $X\bX := X \cup \bX$.
For  $X \in \X$ and $\La \in \B^2$ let  
\[
\begin{array}{ll}
X_{\La} := X \cap \La 
&\text{(restriction of $X$ to $\La$)}, \\
\X_{\La} := \{ X \in \X: X \subset \La \} 
&\text{(set of all configurations in $\La$)  and }\\
N_{\La}(X) := \# X_{\La}  
&\text{(number of particles of $X$ in $\La$).} 
\end{array}
\]
The counting variables $(N_{\La})_{\La \in \B^2}$ generate a $\si$-algebra 
on $\X$, which will be denoted by  $\F_{\X}$. For $\La \in \B^2$ let 
$\F'_{\X,\La}$ be the $\si$-algebra on $\X_{\La}$ obtained by restricting
$\F_{\X}$  to $\X_{\La}$, and let $\F_{\X,\La} := e_{\La}^{-1} \F'_{\X,\La}$ be the 
$\si$-algebra on $\X$ obtained from $\F'_{\X,\La}$ by the restriction mapping 
$e_{\La}: \X \to \X_{\La}, X \mapsto X_{\La}$. The tail $\si$-algebra or
$\si$-algebra of the events far from the origin is defined by 
\[
\F_{\X,\infty} \, := \, \bigcap_{n \ge 1} \F_{\X,\Lan^c}.
\]
Let $\nu$ be the distribution of the Poisson point process on 
$(\X,\F_{\X})$, i.e. 
\[
\int \nu(dX) f(X)  \, = \, e^{- \leb(\La)} \, \sum_{ k \ge 0} \,  
 \frac{1}{k!} \, \int_{{\La}^k} dx_1 \ldots dx_k \, 
 \,  f(\{x_i : 1 \le i  \le k\} ), 
\]
for any $\F_{\X,\La}$-measurable function $f: \X \to \R_+$, where 
$\La \in \Bo_b^2$.  For $\La \in  \Bo^2_b$ and $\bX \in \X$ 
let $\nu_{\La}(.|\bX)$ be the distribution of the Poisson point 
process in $\La$ with boundary condition $\bX$, i.e.
\begin{equation*} 
\int \nu_{\La}(dX|\bX) f(X) \;
= \; \int \nu (dX) f(X_{\La}\bX_{\La^c})
\end{equation*}
for any $\F_{\X}$-measurable function $f: \X \to \R_+$.
It is easy to see that $\nu_{\La}$ is a stochastic kernel 
from $(\X,\F_{\X,{\La}^c})$ to  $(\X,\F_{\X})$.\\

The configuration space of bonds $\E$ is defined to be the set of  all locally 
finite bond sets, i.e.
\[
\E \, := \,  \{B \subset E(\R^2): \, 
\#\{ xx' \in B : xx' \subset \La \} < \infty \, 
\text{ for all } \La \in \B_b^2  \}.
\]
On $\E$ the $\si$-algebra $\F_{\E}$ is defined to be generated by the 
counting variables 
$N_{E}: \E \to \N,  B \mapsto \#(E \cap B)$ $\;(E \in \F_{E(\R^2)})$.\\  
For a countable set $E \in \E$ one can also consider the 
Bernoulli-$\sigma$-algebra $\Bo_E$ on $\E_{E} := \p(E) \subset \E$,  
which is defined to be generated by the family of sets 
$(\{B \subset E: e \in B\})_{e \in E}$. Given a family 
$(p_e)_{e \in E}$ of reals in $[0,1]$ the Bernoulli
measure on $(\E_E,\B_E)$ is defined as the unique probability measure  
for which the events $( \{B \subset E : e \in B \} )_{e \in E}$ 
are independent with probabilities $(p_e)_{e \in E}$.
It is easy to check that the inclusion $(\E_E,\B_E) \to (\E, \F_{\E})$ is 
measurable. Thus any probability measure on  
$(\E_E,\B_E)$ can trivially be extended to $(\E, \F_{\E})$.

%=======================================================================

\subsection{Gibbs measures} \label{secgibbs}

Let $U:\R^2 \to \baR$ be a potential and $z >0$ an 
activity parameter. For finite configurations $X,X' \in \X$ 
we consider the energy terms 
\[
H^U (X) \, := \sum_{x_1x_2 \in E(X)} U(x_1-x_2) \quad \text{ and } \quad  
W^U (X,X') \, := \sum_{x_1 \in X} \sum_{x_2 \in X'} U(x_1-x_2).
\]
The last definition can be extended to infinite configurations  $X'$ 
whenever  $W^U (X,X'_{\La})$ converges as  $\La \uparrow \R^2$ 
through the net  $\B^2_b$. The Hamiltonian of a configuration $X \in \X$ in 
$\La \in \B^2_b$ is given by 
\[ 
H^U_{\La}(X) \, := \, 
H^U (X_{\La}) +  W^U(X_{\La},X_{\La^c}) \, =  
\sum_{x_1x_2 \in E_{\La} (X)}  U(x_1-x_2), 
\] 
where 
\[
E_\La (X) := \{ x_1x_2 \in E(X):  x_1 x_2 \cap \La \ne \emptyset\}.
\]
The integral 
\[ 
Z^{U,z}_{\La}(\bX) \, 
:= \, \int \nu_{\La}(dX|\bX) \, e^{-H^U_{\La}(X)}z^{\#X_\La}
\]
is called the partition function in  $\La \in \B^2_b$ for the boundary condition 
$\bX_{\La^c} \in \X$. In order to ensure that the above objects are 
well defined and the partition function is finite and positive 
we need the following definition: 
\begin{defi} \label{defzul}
A triple $(U,z,\XX)$ consisting of a potential $U:\R^2 \to \baR$, 
an activity parameter $z >0$ and a set of boundary conditions  
$\XX \in \F_{\X,\infty}$ is called admissible if for all 
$\bX \in \XX$ and $\La \in \B^2_b$ the following holds: 
$W^U(\bX_{\La},\bX_{\La^c})$ has a well defined value in $\baR$ and 
 $Z^{U,z}_{\La}(\bX)$ is finite.
\end{defi}
If $(U,z,\XX)$ is admissible, $\La \in \B^2_b$ and $\bX \in \XX$ 
then $W^U(X_{\La},\bX_{\La^c}) \in \baR$ is well defined for every 
$X \in \X$, because $\XX \in \F_{\X,\infty}$ implies $X_{\La}\bX_{\La^c} \in \XX$. 
As a consequence the partition function $Z^{U,z}_{\La}(\bX)$ is well defined. 
Furthermore by definition it is finite and by considering the empty configuration 
one can show that it is positive. The conditional Gibbs distribution 
$\gmu^{U,z}_{\La}(.|\bX)$ in $\La \in \B^2_b$ with boundary condition $\bX \in \XX$ 
is thus well defined by 
\[
\gmu^{U,z}_{\La}(A|\bX) \, := \, 
\frac 1 {Z^{U,z}_{\La}(\bX)}  \int \nu_{\La}(dX|\bX) \, 
e^{-H^U_{\La}(X)} z^{ \#X_\La} 1_A(X) \quad 
\text{ for } \quad  A \in \F_{\X}.
\]
$\gmu^{U,z}_{\La}$ is a probability kernel from $(\XX,\F_{\XX,\La^c})$ to
$(\X,\F_{\X})$. Let 
\[ 
\begin{split}
\G_{\XX}(U,z) \, := \, \{ \mu \in &\prob(\X,\F_{\X}) : \, \mu(\XX) = 1 
\quad \text{ and } \\
&\mu(A|\F_{\X,\La^c}) = 
 \gmu^{U,z}_{\La}(A|.) \text{ $\mu$-a.s. } \alle A \in \F_{\X},
\La \in \B^2_b \}
\end{split}
\]
be the set of all Gibbs measures corresponding to $U$ and $z$ with whole 
weight on boundary conditions in $\XX$. It is easy to see that for any 
probability measure $\mu \in \prob(\X,\F_{\X})$ such that $\mu(\XX)=1$ 
we have the equivalence 
\[ 
\mu \in \G_{\XX}(U,z) \quad \gdw \quad (\mu \otimes \gmu^{U,z}_{\La} = \mu
 \alle \La \in \B^2_b).
\] 
So for every  $\mu \in \G_{\XX}(U,z)$, $f: \X \to \R_+$ measurable 
and $\La \in \B^2_b$ we have 
\begin{equation} \label{gibbsaequ}
\int \mu(dX) \, f(X)\, = \, 
\int \mu(d\bX) \int \gmu_{\La}^{U,z} (dX|\bX) \, f(X). 
\end{equation}
If we consider a fixed potential and a fixed activity we will omit
the dependence on $U$ and $z$  in the notations  $\gmu_{\La}^{U,z}$ and 
$Z_{\La}^{U,z}$.
As a consequence of  \eqref{gibbsaequ} the  hard core $\KU$ of a potential 
$U$ implies that particles are not allowed to get too close to each other, i.e.
for admissible $(U,z,\XX)$ and $\mu \in \G_{\XX}(U,z)$ we have
\begin{equation} \label{hardcore}
\mu(\{X \in \X: \exists x,x' \in X: x\neq x', x-x' \in \KU \})
\, = \, 0. 
\end{equation}
For admissible $(U,z,\XX)$ and a Gibbs measure $\mu \in \G_{\XX}(U,z)$ 
we define the correlation function  $\rho^{U,\mu}$ by 
\[
\rho^{U,\mu}(X) \, = \, e^{-H^U(X)} \int \mu (d\bar{X}) \,  
     e^{- W^U (X,\bar{X})}
\]
for any finite configuration $X \in \X$. If there is a $\xi = \xi(U,z,\XX) \ge 0$ 
such that
\begin{equation} \label{ruellebd}
\rho^{U,\mu}(X) \, \le \, \xi^{\#X} \quad \text{ for all finite }
X \in \X \text{ and all } \mu \in \G_{\XX}(U,z),  
\end{equation}
then we call $\xi$ a Ruelle bound for  $(U,z,\XX)$. Actually we need 
this bound on the correlation function in the following way:
\begin{lem} \label{lekorab}
Let $(U,z,\XX)$ be admissible with Ruelle bound $\xi$.
For every Gibbs measure $\mu \in \G_{\XX}(U,z)$ and every measurable
$f: (\R^2)^m \to \R_+$, $m \in \N$ we have 
\begin{equation} \label{korab}
\int \mu(dX) \sideset{}{^{\neq}}\sum_{x_1,\ldots,x_m \in X} 
f(x_1,\ldots,x_m) \,
\le \, (z\xi)^m  \int dx_1 \ldots dx_m  \, f(x_1,\ldots,x_m).
\end{equation}
\end{lem}
We use $\Sigma^{\neq}$ as a shorthand notation for a multiple sum such that
the summation indices are assumed to be pairwise distinct.

%======================================================================

\subsection{Superstability and admissibility} \label{secadm} 

Now we will discuss some conditions on potentials which imply that $(U,z,\XX)$ 
is admissible and has a Ruelle bound whenever the set of boundary conditions
$\XX$ is suitably chosen. Apart from purely repulsive potentials such as 
the pure hard-core potential considered in Theorem \ref{purehardcore}
we also want to consider superstable potentials in the sense of 
Ruelle, see \cite{R}. Therefore let  
\[
\Gar \, :=  \, r + [-1/2 ,1/2[^2 \subset \R^2 
\qquad (r \in \Z^2)
\] 
be the unit square centred at  $r$ and let 
\[
\Z^2(X) \, := \, \{ r \in \Z^2: N_{\Gar}(X) > 0 \}
\]
be the minimal set of lattice points such that the corresponding squares cover
the configuration. A potential $U:\R^2 \to \baR$ 
is called superstable if there are real constants  $a > 0$ and $b \ge 0$ 
such that for all finite configurations  $X \in \X$  
\[ 
H^U (X) \, \ge \,\sum_{r \in \Z^2(X)} \,[a N_{\Gar}(X)^2 - b N_{\Gar}(X)].
\]
$U$ is called lower regular if there is a decreasing function $\Psi: \N \to \R_+$ 
with $\sum \limits_{r \in \Z^2} \Psi(|r|) <  \infty$ such that
\[
W^U (X,X') \, \ge \, -\sum_{r \in \Z^2(X)} 
\sum_{s \in \Z^2(X')} \, \Psi(|r-s|) \, 
[\frac{1}{2} N_{\Gar}(X)^2 + \frac{1}{2} N_{\Gas}(X')^2] 
\]
for all finite configurations $X,X' \in \X$. So superstability and 
lower regularity give lower bounds on energies 
in terms of particle densities. In order to be able to control these 
densities, a configuration  $X \in \X$ is defined to be tempered if 
\[
\bs(X) \, := \, \sup_{n \in \N} s_n(X) \, < \, \infty, \; \text{ where } \; 
s_n(X) \, := \, \frac{1}{(2n+1)^2} \sum_{r \in \Z^2 \cap \La_{n+1/2}} N_{\Gar}^2(X).
\]
By $\X_t$ we denote the set of all tempered configurations. We note that 
$\X_t \in \F_{\X,\infty}$.  
\begin{lem} \label{lemsup}
Let $z > 0$ and $U:\R^2 \to \baR$ be a translation-invariant pair potential.
\begin{enumerate}
\item[(a)]
If $U$ is purely repulsive, i.e. $U \ge 0$, then  $(U,z,\X)$ is admissible 
with Ruelle bound  $\xi := 1$.
\item[(b)]
If $U$ is superstable and lower regular then $(U,z,\X_t)$ is admissible and
admits a Ruelle bound. 
\end{enumerate}
\end{lem}
The first assertion is a straightforward consequence of the fact that
all energy terms are nonnegative. For the second assertion see \cite{R}.

%======================================================================

\subsection{Conservation of translational symmetry} \label{secsym}

Every $\tauh \in \R^2$ gives a translation on the configuration space $\X$ via  
\[
g_{\tauh}(X) \, := X + \tauh \, := \, \, \{x+ \tauh : x \in X\}.
\]
We say that a measure $\mu$ on $(\X,\F_{\X})$ is 
$\tauh$-invariant if $\mu \circ g_{\tauh}^{-1} = \mu$, 
and $\mu$ is translation-invariant if it is $\tauh$-invariant 
for every $\tauh \in \R^2$. The following lemma gives a sufficient 
condition for the conservation of $\tauh$-symmetry. 
\begin{lem} \label{lemgeorgiikrit}
Let  $(U,z,\XX)$ be admissible, where $U:\R^2 \to \baR$ is a 
translation-invariant potential. If for all cylinder events 
$D \in \F_{\X,\La_{m}}$  $ (m \in \N)$ 
and all Gibbs measures $\mu \in \G_{\XX}(U,z)$ we have  
\begin{equation} \label{georgiikrit}
\mu(D + \tauh) \, + \, \mu(D - \tauh) \, \ge \, \mu(D), 
\end{equation}
then every Gibbs measure  $\mu \in \G_{\XX}(U,z)$ is $\tauh$-invariant.
\end{lem}
We further note that $\R^2$ is generated by the set 
$\{\tau_i e_i:  0 \le \tau_i < 1/2, i \in \{1,2\} \}$, so we only have to 
consider translations of this special form in order to establish 
translation-invariance of a set of Gibbs measures.

%=========================================================================

\subsection{Concerning measurability} \label{secmeas}

We will consider various types of random objects, all of which have
to be shown to be measurable with respect to the considered $\si$-algebras. 
However we will not prove measurability of every such object in detail.
Instead  we will now give a list of operations that preserve
measurability.
\begin{lem} \label{lemeas}
Let $X,X' \in \X$, $B,B' \in \E$,
$x \in \R^2$ and  $p \in \Omega$ be variables, 
where $(\Omega,\F)$ is a measurable space. 
Let $f: \Omega \times \R^2 \to \R$ and 
$g: \Omega \times E(\R^2) \to \R$ be measurable. 
Then the following functions of the given variables are measurable 
with respect to the considered $\si$-algebras: 
\begin{align}
&\sum_{x' \in X} f(p,x'), \quad X \cap X', \quad X \cup X', 
 \quad X \weg X', \quad X+x,\label{measX}\\
&\sum_{b' \in B} g(p,b'), \quad B \cap B', \quad B \cup B', 
 \quad B \weg B',  \quad  B+x, \label{measB}\\
& \inf_{x' \in X} f(p,x'), \quad 
\{x' \in X: f(p,x') = 0\}, \quad C_{X,B}(x),\quad E(X), \label{measinf}\\
&\text{the number of different clusters of $(X,B)$}\label{measncl}.
\end{align}
\end{lem}
Using this lemma and well known theorems, such as the measurability part of 
Fubini's theorem, we can check the measurability of 
all objects considered.

%========================================================================
%
%========================================================================

\section{Proof of the lemmas from Sections \ref{secresult} and 
\ref{secsetting}} \label{secleset}

\subsection{Smoothly approximable potentials: Lemma \ref{lepotchar}}

Let $(U,z,\XX)$, $\xi$, $\psi$, $\tilde{K}$ and $\tU$ (in case (b)) 
be as in the formulation of Lemma~\ref{lepotchar}. By compactness 
of $\KU$ we can choose an $\ep > 0$ such that the $\ep$-enlargement 
$\Ko := (\KU)_{\ep}$ of the hard core $\KU$ has the property
\[
c \, :=  \, 1/(z\xi) - \leb(\Ko \weg \KU) > \, 0.
\]
In case (a) let $U_1 := U$ and in case (b) let $U_1 :=  \tU$.
Let $R \ge 1$ such that 
\[
\Ko \cup \tilde{K} \subset \La_{R} \quad \text{ and furthermore } \quad
 \int_{\La_R^c}  2 \tU(x) |x|^2 dx  \, < \, \frac{c}{2} \quad \text{ in case (b).}
\]
In both cases  $U_1$ serves as an approximation of $U$ on 
$\La_R^c$. Let $C := \overline{\La_{R+1} \weg \Ko}$, $\delta >0$ and 
$f_{\delta}: \R \to \R_+$ be a symmetric smooth probability 
density with support in the $|.|_2$-disc $B_2(\delta)$, e.g. 
\[
f_{\delta}(x) := \frac 1{c_{\delta}} 1_{B_2(\delta)}(x) e^{-(1-|x|_2^2/\delta^2)^{-1}},
\quad \text{ where } c_{\delta} := 
\int_{B_2(\delta)}  e^{-(1-|x|_2^2/\delta^2)^{-1}} dx.
\]
Then
\[
U_2(x) \, :=  \, U \ast f_{\delta} (x) \,
:= \, \int dx' \,  f_{\delta}(x')  U(x-x') 
\]
is a smooth approximation of $U$ on $C$. By continuity of $U$ and compactness 
of $C$ a sufficiently small $\delta$ guarantees
\[
|U_2(x) - U(x)| \, <  \, c'  := \, \frac{c}{4\leb(C)} 
\quad \text{ for  } x \in C.
\]
Let  $g: \R^2 \to [0,1]$ be a smooth symmetric function such that 
$g = 0$ on $\La_R$ and $g = 1$ on $\La_{R+1}^c$. Now we can define 
$\bU, u: \Ko^c \to \R$ by  
\[
\bU \, := \, (1 - g) (U_2 + c') + g U_1 \quad \text{ and } \quad 
u \, := \, \bU - U. 
\]
It is easy to verify that the constructed objects have all the properties 
described in Definition  \ref{defGapprox} in both cases (a) and (b).

%=====================================================================

\subsection{Property of the Ruelle bound: Lemma \ref{lekorab}} 

For every $n \in \N$, every measurable $g: \X_{\La_n} \to \R_+$ 
and every $\bX \in \XX$ we have 
\begin{displaymath}
\begin{split}
\int &\nu_{\Lan}(dX|\bX)  \sideset{}{^{\neq}}\sum_{x_1,\ldots,x_m \in 
    X_{\La_n}} f(x_1,\ldots,x_m) \, g(X) \\ 
&= \, \int_{{\Lan}^m} dx_1 \ldots dx_m 
  \, f(x_1,\ldots,x_m) \,
  \int \nu_{\Lan}(dX'|\bX)\, g( \{x_1, \ldots, x_m \} X').
\end{split}
\end{displaymath}
Combining this with \eqref{gibbsaequ}, the definition of the conditional 
Gibbs distribution and the definition of the correlation function we get 
\begin{displaymath}
\begin{split}
&\int \mu(dX) \sideset{}{^{\neq}}\sum_{x_1,\ldots, x_m \in X_{\Lan}} 
 f(x_1,\ldots ,x_m)\\
&= \,  \int \mu(d\bX) \, \frac{1}{Z^{U,z}_{\La_{n}}(\bX)}
   \, \int \nu_{\Lan} (dX|\bX) \sideset{}{^{\neq}} 
   \sum_{x_1,\ldots ,x_m \in X_{\Lan}}  f(x_1,\ldots ,x_m) \, 
    e^{- H^U_{\La_{n}}(X)}z^{ \# X_{\Lan}}\\   
&= \, \int_{{\Lan}^m} dx_1 \ldots  dx_m \, f(x_1,\ldots ,x_m)\, z^m 
   \, \rho^{U,\mu} (\{x_1,\ldots ,x_m\}).  
\end{split} 
\end{displaymath} 
Now we use  \eqref{ruellebd} to estimate the correlation function by 
the Ruelle bound  $\xi$.  Letting  $n \to \infty$ the assertion follows
from the monotone limit theorem. 

%=======================================================================

\subsection{Sufficient condition: Lemma  \ref{lemgeorgiikrit}}

The lemma can be shown exactly as Proposition (9.1) in \cite{G} 
and we will only outline the proof:
We first note that $(\X,\F_{\X})$ is a standard Borel space, 
which  follows from \cite{DV}, Theorem A2.6.III. Hence the
point particle version of Theorem (7.26) in \cite{G} implies 
that every Gibbs measure can be decomposed into extremal Gibbs measures. 
Thus without loss of generality we may assume $\mu$ to be extremal. 
Suppose now that $\mu$ is not $\tauh$-invariant, i.e. 
$\mu \circ g_{\tauh}^{-1} \neq \mu$, which also implies 
$\mu \circ g_{\tauh} \neq \mu$. As the extremality of $\mu$ implies 
the extremality of  $\mu \circ g_{\tauh}^{-1}$ and $\mu \circ g_{\tauh}$, 
the point particle version of Theorem (7.7) guarantees the existence 
of sets $A_-,A_+ \in \F_{\X,\infty}$ such that  $\mu \circ g_{\tauh}^{-1}(A_-)=0$, 
$\mu \circ g_{\tauh}(A_+)=0$ and $\mu(A_-)=\mu(A_+)= 1$. Hence for 
$A := A_- \cap A_+$ we have 
\[
\mu \circ g_{\tauh}(A) + \mu \circ g_{\tauh}^{-1}(A) \, =\,0 \,< \,1 \,= \,\mu(A).
\]
On the other hand by assumption \eqref{georgiikrit} we know that 
$\mu \circ g_{\tauh} + \mu \circ g_{\tauh}^{-1} \ge \mu$ on the algebra of 
all cylinder events. By the monotone class theorem this inequality 
even holds on all of $\F_{\X}$, 
which contradicts the above inequality.

%=======================================================================

\subsection{Measurability: Lemma \ref{lemeas}}

Details concerning measurability of functions of point processes 
can be found in \cite{DV}, \cite{K} or \cite{MKM}, for example . 
The first part of \eqref{measX} is the 
measurability part of Campbell's theorem.
For the rest of \eqref{measX} it suffices to observe that for 
$\La \in \Bo^2_b$ we have  
\[ \begin{split}
&N_{\La}(X \cap X') 
 = \sum_{x \in X} \sum_{x' \in X'} 1_{\{x = x' \in \La\}}, \; 
 N_\La(X \weg X') =  N_\La(X) -  N_\La(X \cap X'),\\
&N_\La(X+x) = \sum \limits_{x' \in X} 1_\La(x'+x) \; \text{ and } \;   
N_\La(X \cup X') =  N_\La(X) +  N_\La(X' \weg X). 
\end{split} \]
\eqref{measB} can be proved similarly.
For $c \in \R$, $\La \in \B^2_b$, $x' \in \R^2$ and   $L \in \F_{E(\R^2)}$
\[ \begin{split} 
&\inf_{x' \in X} f(p,x') < c \, 
\gdw \, \sum_{x' \in X} 1_{\{f(p,x')<c\}} \ge 1 , \\  
&N_\La(\{x' \in X: f(p,x') = 0\}) \, 
 =  \, \sum_{x' \in X} 1_{\{f(p,x')= 0, x' \in \La\}}, \\  
&N_\La(C_{X,B}(x)) \, = \, \sum_{x' \in X} 
 1_{\{x' \in C_{X,B}(x), x' \in \La\}},\\
&x' \in C_{X,B}(x) \, \gdw \,  \sum_{m \ge 0} \sum_{x_0,\ldots ,x_m \in X} 
 1_{\{x=x_0,x'=x_m\}}\prod_{i=1}^m 1_{\{ x_i x_{i+1} \in B  \}} \ge 1 \quad \text{ and }\\
&N_L(E(X)) \, = \, \frac 1 2 \sum_{x_1 \in X} \sum_{x_2 \in X \weg \{x_1\}}
 1_{L}(x_1x_2).  
\end{split} \]
Using these relations, the measurability of the terms in \eqref{measinf} follows 
easily. For \eqref{measncl} it suffices to observe that there are at most 
$k$ different clusters of $(X,B)$ iff  
\[
\sum_{x_1,\ldots ,x_k \in X} 1_{\{ X \weg (C_{X,B}(x_1) \cup \ldots  
 \cup  C_{X,B}(x_k)) = \emptyset \}} \ge 1. 
\]

%========================================================================
%
%========================================================================

\section{Proof of Theorem \ref{purehardcore}: Main steps} 
\label{secproofcore}

\subsection{Basic constants} \label{constantshc}

Let $z>0$. Let $|.|_h$ be a norm on $\R^2$ and $U := U_{hc}$ the corresponding 
pure hard-core potential. As $U$ is purely repulsive 
we know that $(U,z,\X)$ is admissible with Ruelle bound $\xi := 1$
by Lemma \ref{lemsup}, part (a). Let $K := \KU$ and 
$\ep >0$. If we choose $\ep$ sufficiently small we have 
\begin{equation} \label{dechc}
\ex \, := \, \leb(\Kep \weg \KU) \, < \, \frac{1}{z\xi}, 
\end{equation}
where $\Kep$ is the $\ep$-enlargement of $\Ko$. Let $\fc:\R^2 \to \R$  
be a function such that 
\[
\text{ $\fc$ is smooth, $\quad \fc = 0$ on $\Ko \quad $ and 
$\quad \fc = 1$ on $(\Kep)^c$.}
\] 
Furthermore we need the following finite constants: 
\begin{equation} \label{bpsihc}  
\cK \,  := \, \sup\{|x|: x  \in \Kep \} \quad \text{ and }\quad 
\cf \,  := \, \sup\{|\fc'(x)|: x \in \R^2  \}. 
 \end{equation}
On $\R^2$ let $\le$ be the lexicographic order and let the partial order 
$\lear$ be defined by 
\[
 (r_1,r_2) \, \lear \, (r_1',r_2') \quad 
:\Leftrightarrow \quad r_1 \le r_1', r_2 = r_2'.
\]
In order to show the conservation of translational symmetry we 
fix a Gibbs measure  $\mu \in \G_{\XX}(U,z)$ and  
a cylinder event $D \in \F_{\X,\La_{n'-1}}$ where $n'\in \N$, 
see Subsection \ref{secsym}. As mentioned there it suffices
to consider translations $\tau e$, where $\tau  \in [0,1/2]$ 
and $e = e_1$ or $e_2$. Hence we fix $\tau  \in [0,1/2]$, and  by symmetry 
we may assume that $e = e_1$. We also fix an arbitrarily small real 
$\de > 0$ in order to control probabilities close to $0$. 
As all the above objects are fixed for the whole proof we will ignore 
dependence on them in our notations.

%=========================================================================

\subsection{Generalised translation} \label{tgthc} 

\begin{figure}[!htb] 
\begin{center}
\psfrag{1}{$0$}
\psfrag{2}{$R$}
\psfrag{3}{$n$}
\psfrag{4}{$s$}
\psfrag{5}{$\tau$}
\psfrag{6}{$\tn(s)$}
\includegraphics[scale=0.4]{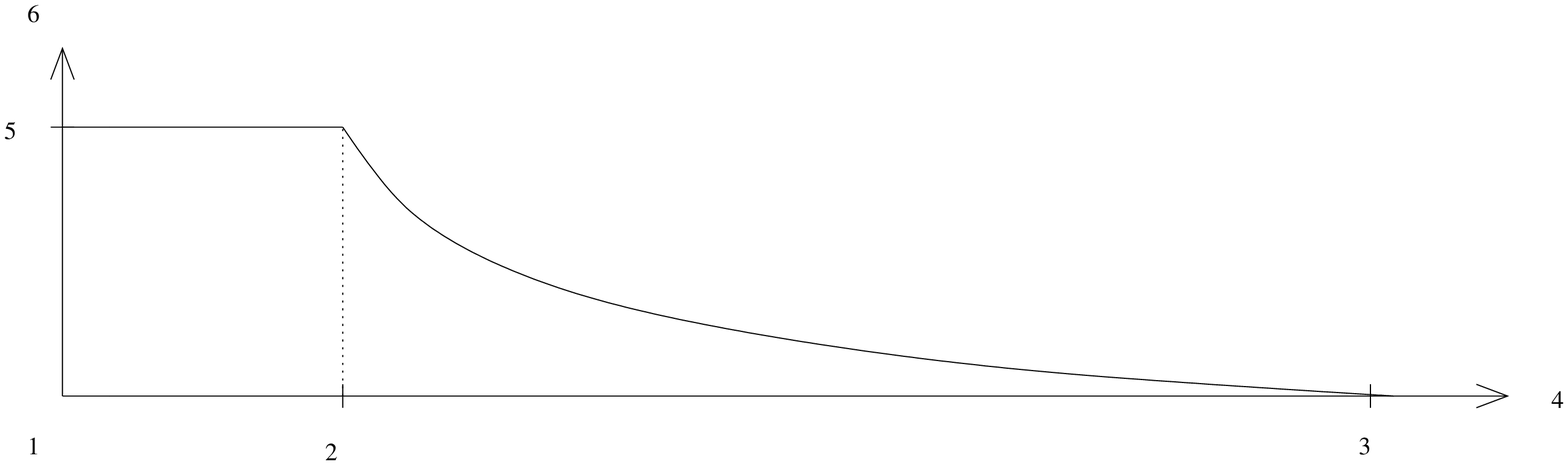}
\end{center}
\caption{Graph of $\tn$}\label{figtauhc}
\end{figure}

Let $n > n'$ and  $X \in \X$. We consider the bond set 
\[
\Kepn \, := \, 
\{ x_1x_2 \in E(X): x_1x_2 \cap \Lan \neq \emptyset, x_1-x_2 \in \Kep \}.
\]
Every time we use this notation it will be clear from the
context which configuration $X$ it refers to. Note that $\Kepn$ is finite
as $X$ is locally finite and $\Kep$ is bounded. For a bounded set 
$\La \in \Bo^2_b$ let
\begin{displaymath} 
\rnx(\La) \, = \, \sup\{ |y'| : y' \in \hcny(\La) \} 
\end{displaymath}
denote the range of the corresponding $\Kepn$-cluster.
In the following lemma we consider the case $\La = \La_{n'}$, 
where $n' \in \N$ is the number fixed in Section \ref{constantshc}. 
\begin{lem} \label{lereichhc}
We have $ \quad \sup \limits_{n > n'} \int \mu (dX) 
   \, \rnx(\La_{n'}) \,  < \, \infty$.
\end{lem}
By the Chebyshev inequality we therefore can choose an 
integer $R > n'$, such that for every $n >n'$ 
we have
\[
\mu  (\Gna) \, \ge \, 1- \frac \de 2 \quad \text{ for } \quad 
\Gna \, := \, \{X \in \X: \rnx(\La_{n'}) <  R\} \, 
\in \, \F_{\X}.
\]
For $n > R$ we define the functions
\[
q: \R_+ \to \R, \quad Q: \R_+ \to \R, \quad r: \R \times \R_+ \to \R \quad 
\text{ and } \quad  \tn: \R \to \R  \quad \text{ by} 
\]
\begin{displaymath}
\begin{split}
q(s) \, &:= \, \frac{1}{1 \vee (s \log(s))},
\hspace{ 1,3 cm}Q(k) \, := \, \int_0^k q(s) ds, \\
r(s,k) \, &:= \, \int_{(s \vee 0) \wedge k} ^k \frac{q(s')}{Q(k)} ds', \qquad  
\tn(s) \, := \,  \tau \, r(s-R,n-R).
\end{split}
\end{displaymath}
For a sketch of the graph of $\tn$ see Figure \ref{figtauhc}. 

\noindent
Some important properties of $\tn$ are the following:
\begin{equation} \label{tauinnenaussenhc}
\begin{split}
\tn(s) = \tau \, &\text{ for } \, s \le R,  \quad 
\tn(s) = 0 \, \text{ for } \, s \ge n\\
&\text{ and $\tn$ is decreasing.}
\end{split}
\end{equation}
For $X \in \X$ and $x \in X$ we define $\anx(x)$ to be the point of 
$\hcny(x)$ with maximal $|.|$-distance to the origin. 
(If there is more than one such point we choose the maximal one 
with respect to the lexicographic order for the sake of definiteness.)  
Then \eqref{tauinnenaussenhc} implies 
\[
|\anx(x)| \; \ge \; |x| \quad \text{ and } \quad 
\tn(|\anx(x)|) \; 
= \; \min \{ \tn(|x'|) : x' \in \hcny(x) \}.
\]
The transformation $\eTn^0 : \R^2 \to \R^2$, $\eTn^0(x) := x+\tn(|x|)\einh$
can also be viewed as a transformation on $\X$, such that 
every point $x$ of a configuration $X$ is translated the distance 
$\tn(|x|)$ in direction $\einh$. We would like to use this generalised 
translation $\eTn^0$ as a tool for our proof just as in \cite{FP1} 
and \cite{FP2}. 

%=========================================================

\subsection{Good configurations}

In order to deal with the hard core we will replace the above translation 
$\eTn^0$ by a transformation  
\[
\Tn: \X \to \X 
\]
which is required to have the following properties: 
\begin{enumerate}
\item [(1)] For $X \in \X$ the transformed configuration  $\tX = \Tn(X)$  
is constructed by translating every $x \in X$ a certain distance  
$\tnx(x)$ in direction  $\einh$. We note that  we do not require 
the particles to be translated independently. 
\item [(2)] Particles in the inner region $\La_{n'-1}$ are translated by 
$\tau \einh$,
and particles in the outer region  ${\Lan}^c$ are not translated at all.
\item [(3)] $\Tn$ is bijective, the density of the transformed process with 
respect to the untransformed process under the measure $\nu$ 
can be calculated explicitly and we have a  suitable estimate on this density. 
\item [(4)] The Hamiltonian $H^{U}_{\Lan}(X)$ is invariant under $\Tn$, 
i.e. particles within hard core distance remain within hard core distance
and particles at larger distance remain at larger distance. 
\end{enumerate} 
Property (2) implies that the translation of the chosen cylinder event $D$ 
is the same as the transformation of $D$ by $\Tn$. Properties (3) and (4) 
imply that the density of the transformed process with 
respect to the untransformed process under the measure $\mu$
can be estimated. Therefore a transformation with these properties 
seems to be a good tool for proving
\eqref{georgiikrit}. However, in general it is difficult to 
construct a transformation with all the given properties. For example
properties (2) and (4) cannot both be satisfied if $X$ is a configuration 
of densely packed hard-core particles. 
If  $n > R$ and $X \in \Gna$ then such a situation can not occur, and 
by Lemma \ref{lereichhc} this is the case with high probability. 
Similar problems arise for the other properties, so we will content ourselves 
with a transformation satisfying the above properties only for  
configurations $X$ from a set of good configurations 
\begin{equation} \label{goodhc}
 \Gn  \, :=  \, \big\{ X \in \Gna: 
  \sum \limits_{i=1}^3 \Sigma_i(n,X) <   1 \big\} \, 
\in \, \F_{\X}.
\end{equation} 
The functions $\Sigma_i(n,X)$ will be defined whenever we want
good configurations to have a certain property. In Lemma
\ref{lebadsmallhc} we then will prove that the set of good configurations 
$\Gn$ has probability close to $1$ when $n$ is big enough.  
Up to that point we consider a fixed $n \ge R+1$. 

%=====================================================================

\subsection{Modifying the generalised translation} \label{secmodgen}

With a view to properties (1) and the second part of (2) 
we define the transformation $\Tn: \X \to \X$ by  
\[ 
\Tn(X) \, := \, X_{\Lan^c} \cup \{\pnx^k + \taunx^k \einh: 1 \le k \le m(X)\} =
\{x + \tnx(x) \einh : x \in X \}
\]
for every $X \in \X$, where $m(X) := \# X_{\Lan}$, 
$\{\pnx^k: 1 \le k \le m(X)\} = X_{\Lan}$, 
$\taunx^k$ is the translation distance of $\pnx^k$ and the 
translation distance function  $\tnx: X \to \R$ is defined by 
$\tnx(x) := 0$ for $x \in X_{\Lan^c}$ and $\tnx(\pnx^k) := \taunx^{k}$ for 
$1 \le k \le m(X)$. 
\begin{figure}[!htb] 
\begin{center}
\psfrag{1}{$\etaun^0 = 0$}
\psfrag{2}{$\etaun^0 = 0$}
\psfrag{3}{$\etaun^1$}
\psfrag{4}{$\etaun^2$}
\psfrag{5}{$\etaun^3$}
\psfrag{7}{$\etaun^4 = \tau$}
\psfrag{a}{$\Lan$}
\psfrag{b}{$\La_R$}
\includegraphics[scale=0.65]{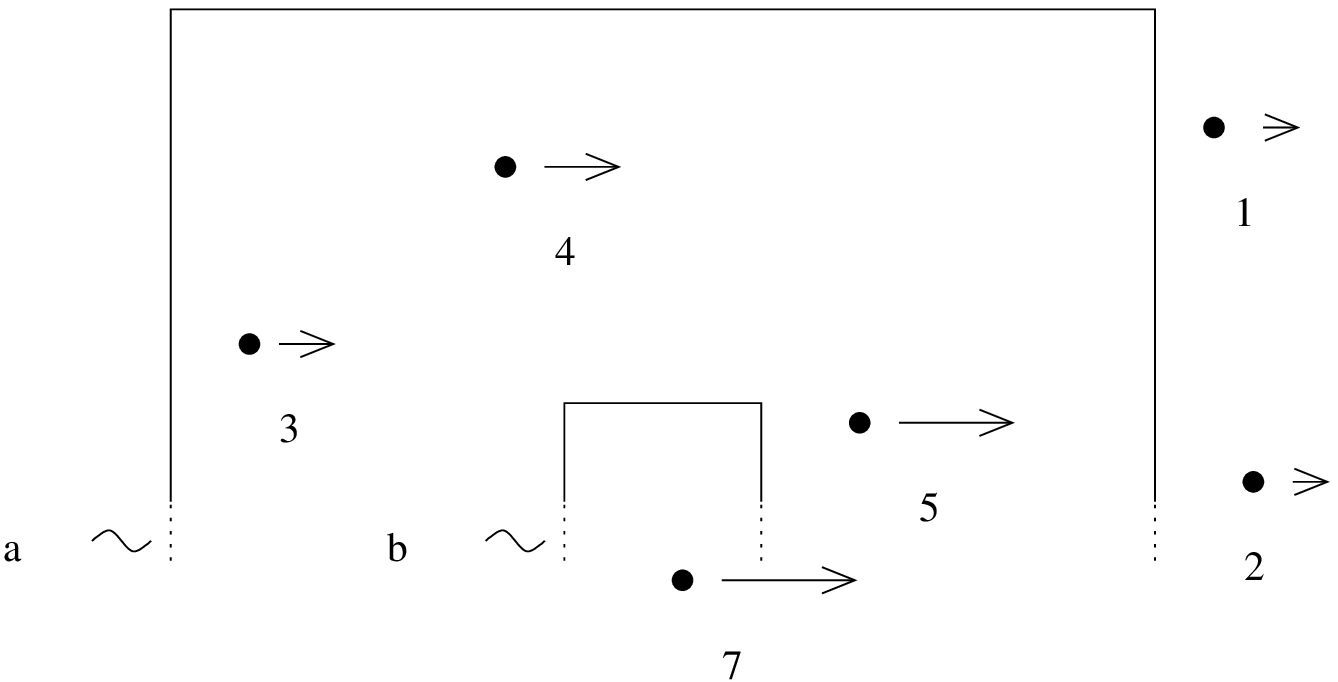}
\end{center} 
\caption{Every point $P^k_n$ is translated by $\etaun^k\einh$} 
\label{figclushc}
\end{figure}
We are left to identify the points $\pnx^k$ of $X$ 
and their translation distances $\taunx^k$. In order to simplify notation 
we will omit the dependence on $X$ in $m(X)$, $\tnx, \pnx^k$ and $\taunx^k$ 
if it is clear which configuration is considered. 
In our construction we would like to ensure that the points $\epn^k$ 
are ordered in a way such that 
\begin{equation} \label{monohc}  
0 =: \etaun^0 \le \etaun^1 \le \ldots \le \etaun^{m}. 
\end{equation}
This relation will be an important tool for showing the bijectivity of the 
transformation as required in property (3) of the last subsection. 
As required in (4) we also would like to have
\begin{align}
&x_1,x_2 \in X,\,  x_1-x_2 \in \Ko  \;
\Rightarrow \; \tnx(x_1) \,= \,\tnx(x_2), \label{verconhc}\\
&x_1,x_2 \in X, \, x_1-x_2 \notin \Ko \; \Rightarrow \; 
(x_1+ \tnx(x_1)\einh)-(x_2+ \tnx(x_2)\einh) \notin \Ko. \label{grepshc} 
\end{align}
With these properties in mind we will now give a recursive definition 
of $\epn^k$ and $\etaun^k$ for a fixed configuration $X \in \X$ 
using a translation distance function $\etn^k := \tnx^k: \R^2 \to \R$ 
in each step. In the k-th construction step  $(1 \le k \le m)$ let 
\[
\begin{split} 
&\etn^{k} \, := \,\etn^{0} \wedge \bigwedge_{0 \le i < k} m_{\epn^{i},\etaun^{i}}\, 
= \, \etn^{k-1} \wedge m_{\epn^{k-1},\etaun^{k-1}},\\
&\text{ where } \quad \etn^0 \, := \, \tn(|.|) \quad  \text{ and } \quad 
m_{\epn^{0},\etaun^{0}} \, := \,  \bigwedge_{x \in X_{\Lan^c}}  m_{x,0}. 
\end{split} \]
The auxiliary functions $\mxt$ will be defined later. 
\begin{figure}[!htb]
\begin{center}
\psfrag{1}{$\La_R$}
\psfrag{2}{$\Lan$}
\psfrag{3}{$\epn^2$}
\psfrag{4}{$\epn^1$}
\psfrag{a}{$\epn^3$}
\psfrag{6}{$\etaun^1$}
\psfrag{5}{$\etaun^2$}
\psfrag{b}{$\etaun^3$}
\psfrag{8}{$m_{\epn^1,\etaun^1}$}
\psfrag{7}{$m_{\epn^2,\etaun^2}$}
\psfrag{c}{$m_{\epn^3,\etaun^3}$}
\psfrag{9}{$\etn^0$}
\psfrag{10}{$\tau$}
\psfrag{11}{$\R^2$}
\psfrag{12}{$\R$}
\includegraphics[scale=0.45]{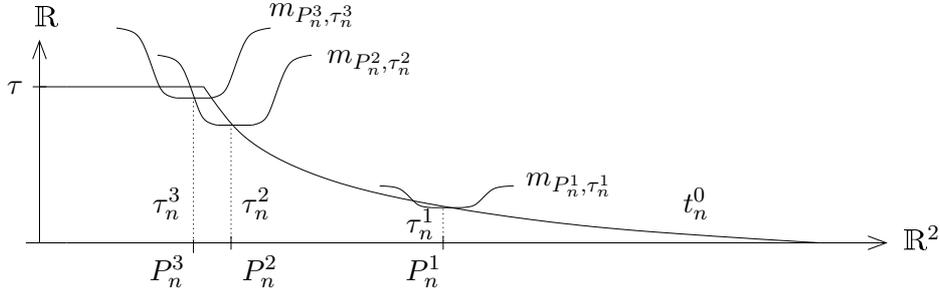}
\end{center}
\caption{Construction of $\etn^k$}  \label{figtrdhc}
\end{figure}

\noindent
Let $\epn^k$ be the point of $X_{\Lan} \weg \{\epn^1, \ldots, \epn^{k-1}\}$ 
at which the minimum of $\etn^{k}$ is attained. 
If there is more than one such point then take the smallest point 
with respect to the lexicographic order for the sake of definiteness. 
Let  $\etaun^k := \etn^{k}(\epn^k)$  be the corresponding minimal value of 
$\etn^k$ and  $\Tnxn^k := \eTn^k := id + \etn^k \einh$.\\
$\etn^k$ is defined to be 
$\etn^0$ modified by local distortions $\mxt$. On the one hand we have thus
ensured  that $\etn^k-\etn^0$ is small, i.e. $\etaun^k \approx \tn(|\epn^k|)$, 
which will give us hold on the density in property (3). On the other hand 
the auxiliary functions of the form $\mxt$ slow down the translation locally  
near every point $x'$ with known translation distance $t$, see Figure \ref{figtrdhc}. 
This will ensure properties \eqref{verconhc} and \eqref{grepshc}. 
\begin{figure}[!htb] 
\begin{center}
\psfrag{1}{$x'+\Ko$}
\psfrag{2}{$x' + \Kep$}
\psfrag{3}{$t$}
\psfrag{4}{$\infty$}
\psfrag{5}{$\hxt$}
\psfrag{6}{$x \in \R^2$}
\psfrag{7}{$\mxt(x)$}
\includegraphics[scale=0.35]{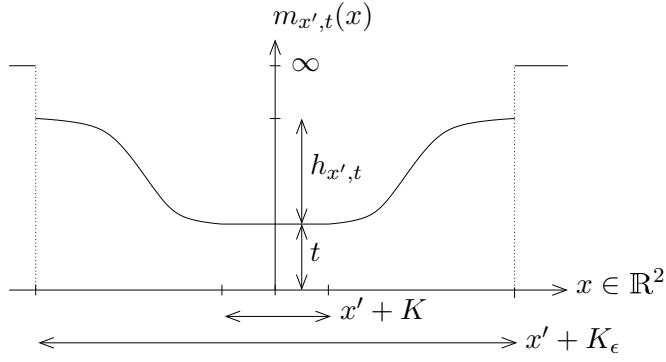}
\end{center}
\caption{One-dimensional sketch of the graph of $\mxt$} \label{figmxthc}
\end{figure}

\noindent 
For $x' \in \R^2$ and $t \in \R$ let the auxiliary function 
$\mxt: \R^2 \to \bar{\R}$ be given by 
\[
\begin{split}
\mxt(x) \, := \, &\Bigg\{ 
\begin{array}{cl}
t & \text{ if } \,  
 \hxt \cf > \frac{1}{2}\\
 t +  \hxt \fc(x-x')  + \infty \, 1_{\{ \fc(x-x')=1\}} \qquad  
& \text{ else},  
\end{array}\\  
&\text{ where } \quad \hxt \, :=  \, |\tn(|x'| - \cK) - t|.
\end{split} 
\] 
Note that the first case in the definition of $\mxt$ has been introduced
in order to bound the slope of $\mxt$. In Section \ref{secmeig} we will 
show important properties of this auxiliary function, but for the moment 
we will content ourselves with the intuition given by Figure \ref{figmxthc}.
Using Lemma \ref{lemeas} one can show that all above objects are measurable 
with respect to the considered $\si$-algebras. In the rest of this section 
we will convince ourselves that the above construction has indeed all the
required properties. 
\begin{lem}\label{leconshc}
The construction satisfies \eqref{monohc}, \eqref{verconhc} and \eqref{grepshc}.
\end{lem}
\begin{lem} \label{leinnenaussenhc}
For good configurations $X \in \Gn$ we have    
\begin{equation}
(\Tn X - \tau \einh)_{\La_{n'-1}} \, = \, X_{\La_{n'-1}}\quad \text{ and } \quad
(\Tn X)_{{\Lan}^c} \, = \, X_{{\Lan}^c}.  \label{innenaussenhc}
\end{equation}
\end{lem}
\begin{lem} \label{lebijhc}
The transformation $\Tn: \X \to \X$ is bijective.
\end{lem}
Actually in the proof of Lemma \ref{lebijhc} we construct the inverse of $\Tn$. 
This is needed in the proof Lemma \ref{ledensphc}, where we will show  for every 
$\bX \in \X$  that  $\nu_{\Lan}(.|\bar{X})$ is absolutely continuous with 
respect to $\nu_{\Lan}(.|\bar{X}) \circ \Tn^{-1}$ with density $\ph \circ \Tn^{-1}$, 
where 
\begin{equation} \label{densdefhc}
\ph (X) \, := \, 
\prod_{k = 1}^{m(X)} \big|1 +  \partial_{1} \tnx^{k} (\pnx^k)\big|. 
\end{equation}
The proof will also show that definition \eqref{densdefhc} makes sense 
$\nu_{\Lan} (\,. \, |\bar{X})$-a.s., in that the considered derivatives exist. 
\begin{lem} \label{ledensphc}
For every  $\bar{X} \in \X$ and every $\F_{\X}$-measurable function $f \ge 0$ 
\begin{equation} \label{denshc}
\int d\nu_{\Lan}(.|\bar{X}) \, (f \circ \Tn \cdot \ph) \,
= \, \int d\nu_{\Lan}(.|\bar{X}) \, f.
\end{equation}
\end{lem}
Considering \eqref{georgiikrit} we also need  the backwards translation. So let
$\iTn$ and $\bph$ be defined analogously to the above objects, 
where now $\einh$ is replaced by $-\einh$. The previous lemmas apply analogously 
to this deformed backwards translation. We note that $\iTn$ is not the inverse 
of $\Tn$. 

%=======================================================================

\subsection{Final steps of the proof}  \label{Potthc}

From \eqref{gibbsaequ} and  Lemma \ref{ledensphc} we deduce  
\[ \begin{split}
\mu & (\Tn (D \cap \Gn))\\ 
&= \,  \int \mu(d\bX) \frac 1 {Z_{\Lan}(\bar{X})} 
 \int \nu_{\Lan} (dX|\bX) \, 1_{\Tn(D \cap \Gn)} (X) \, z^{\# X_{\Lan}} \, 
 e^{-H^{U}_{\Lan}(X)} \\ 
&= \,  \int \mu(d\bX) \frac 1 {Z_{\Lan}(\bar{X})} 
  \int \nu_{\Lan}(dX|\bX)\\ 
&\hspace{2.7 cm} 1_{\Tn(D \cap \Gn)} \circ \Tn(X) \, z^{ \# (\Tn X)_{\Lan}} \,
  e^{-H^{U}_{\Lan}(\Tn X)} \, \ph(X).
\end{split} \]
By  Lemma \ref{lebijhc}  $~\Tn$ is bijective, by \eqref{innenaussenhc}
$~\# (\Tn X)_{\Lan} = \# X_{\Lan}$ and by \eqref{verconhc} and \eqref{grepshc} 
we have $H^{U}_{\La_{n}}(\Tn X) = H^{U}_{\Lan}(X)$. Hence the above integrand 
simplifies to 
\[
 1_{D \cap \Gn} (X)  \, z^{\# X_{\Lan}} \, 
  e^{ - H^{U}_{\Lan}(X)} \ph (X),
\]
and we have an analogous expression for the backwards transformation  $\iTn$. 
So
\[ \begin{split}
\mu  (\iTn &(D \cap \Gn))\, 
 + \,  \mu  (\Tn (D \cap \Gn)) \, 
 - \, \mu (D \cap \Gn) \\ 
&= \,  \int \mu(d\bX) \frac{1}{Z_{\La_{n}}(\bar{X})} 
  \int \nu_{\Lan}(dX|\bX) \, 1_{D \cap \Gn} (X)  
 \, z^{\# X_{\Lan}}   e^{ - H^{U}_{\Lan}(X)} \\  
& \hspace{6.5 cm} \times  \Big[ \bph(X) + \ph(X) - 1 ].
\end{split} \]
We note that for $X \in \Gn$ we have 
\[
\bph(X) + \ph(X) \, \ge \, 2 \, (\bph(X)\ph(X))^{\frac{1}{2}}\,
\ge \, 2 \, e^{- \frac 1 2 } \,
\ge \, 1,  
\]
where we have used the arithmetic-geometric-mean inequality in the first step
and the following estimate in the second step:  
\begin{lem} \label{ngrosshc}
For $X \in \Gn$ we have     
\begin{equation} \label{dichtenkleinhc}
\log \bph(X)  + \log \ph(X) \, \ge \, -1. 
\end{equation}
\end{lem}
Hence we have shown that 
\begin{equation} \label{hingeorgiihc}
\mu (\iTn (D \cap \Gn)) \, + \, 
\mu (\Tn (D \cap \Gn)) \, 
\ge \, \mu  (D \cap \Gn) .
\end{equation}
In \eqref{hingeorgiihc} we would like to replace $D \cap \Gn$ by $D$, 
and for this we need $\Gn$ to have high probability:
\begin{lem} \label{lebadsmallhc}
If $n \ge R + 1$ is chosen big enough, then $\mu  (\Gn^c) \le \de$.
\end{lem}
For the proof of Theorem \ref{purehardcore} we choose such an $n \ge R+1$. 
Because of $D \in \F_{\X,\La_{n'-1}}$ and \eqref{innenaussenhc} we have  
\[
\alle X \in D \cap \Gn: \quad (\Tn X - \tau \einh)_{\La_{n'-1}} \in D, 
\quad \text{ i.e. } \Tn X \in D + \tau \einh,
\]
and an analogous result for the backwards transformation. Hence 
\[
\Tn(D\cap\Gn) \subset D + \tau \einh \quad \text{ and } \quad 
\iTn(D\cap\Gn) \subset D - \tau \einh.
\] 
Using these inclusions and Lemma \ref{lebadsmallhc} we deduce from 
\eqref{hingeorgiihc} 
\[
\mu(D - \tau \einh) \, + \, 
\mu(D + \tau \einh) \, \ge \, \mu(D) \, - \,  \de.
\]
$\de > 0$ was chosen to be an arbitrary positive real, so we get the estimate 
\eqref{georgiikrit} by taking the limit  $\delta \to 0$. 
Now the claim of the theorem follows from Lemma~\ref{lemgeorgiikrit}.

%=========================================================================
%
%=========================================================================

\section{Proof of the lemmas from Section \ref{secproofcore}}
\label{secleproofcore}

\subsection{Cluster bounds: Lemma \ref{lereichhc}} \label{seclereichhc}

For $n>n'$ and $X \in \X$ we want to estimate $\rnx(\La_{n'})$.  
For any path $x_0,...,x_m$ in the graph 
$(X,\Kepn)$ such that $x_0 \in \La_{n'}$ we have 
\[
|x_m| \, \le \,  |x_0| + \sum_{i=1}^m |x_i-x_{i-1}| \,  
\le \, n' + m \cK. 
\] 
By considering all possibilities for such paths we obtain
\[
\rnx (\La_{n'}) \, \le \, n' +  \sum_{m \ge 1}
 \quad \sideset{}{^{\neq}} \sum_{x_0,\ldots ,x_m \in X}   1_{\{x_0 \in \La_{n'} \}}  
 m \cK \prod_{i=1}^m   1_{ \{x_i x_{i-1} \in \Kepn \} }.
\]
Using the hard core property  \eqref{hardcore} and Lemma \ref{lekorab} we get 
\[ \begin{split}
R_n \, &:= \, \int \mu(dX) \rnx (\La_{n'}) \, - \, n'\\
&\le \, \sum_{m \ge 1} \int \mu(dX) \, \sideset{}{^{\neq}}
 \sum_{x_0,\ldots ,x_m \in X} 1_{\{x_0 \in \La_{n'} \}}
  m \cK \prod_{i=1}^m   1_{\Kep \weg \KU}(x_i-x_{i-1})\\
&\le \, \sum_{m \ge 1} \, (z\xi)^{m+1} \, 
   \int dx_0 \ldots  dx_{m}   1_{\{x_0 \in \La_{n'} \}}
  m \cK \prod_{i=1}^m  1_{\Kep \weg \KU} (x_i-x_{i-1}).
\end{split} \] 
By  \eqref{dechc} we
can estimate the integrals over $dx_i$ in the above expression 
beginning with  $i =m$. This gives $m$ times a factor 
$\ex$ and the integration over  $dx_0$ 
gives an additional factor  $\la^2(\La_{n'}) = (2n')^2$. Thus   
\[
  R_n \, \le  \, (2n')^2 z\xi \cK \sum_{m \ge 1}  
  m (\ex z\xi )^{m}  \, < \, \infty,
\]
where the last sum is  finite because $\ex z\xi < 1$.
 
%====================================================================

\subsection{Properties of the auxiliary function } \label{secmeig}

Let  $f: I \to \R$ be a function on an interval. $f$ is called 
 $1/2$-Lipschitz-continuous  if   
\[
|f(r)-f(r')| \, \le \, \frac 1 2  |r-r'| \quad 
\text{ for all  } r,r' \in I.
\]
$f$ is called piecewise continuously differentiable
if it is continuous and if
\[
\gibt \text{ countable and closed }M \subset I: \; \text{ $f$ is continuously 
differentiable on $I \weg M$.} 
\]
As $M$ is closed, the connected components of $I \weg M$ are countably 
many intervals. For a strictly monotone piecewise continuously differentiable 
transformation $f$ on $\R$ we can apply the Lebesgue transformation theorem: 
The derivative $f'$ is well defined $\la^1$-a.s. and for every  
$\B^1$-measurable function $g \ge 0$ we have  
\begin{equation} \label{lebtrsatz}
\int g(f(x)) |f'(x)| dx \, = \, \int g(x') dx'. 
\end{equation}
The above properties are inherited as follows: 
\begin{lem} \label{leanalmin} 
Let $f_1,f_2: I \to \R$ be functions on an interval $I$.
\begin{enumerate}
\item [(a)]If $f_1$ and $f_2$ are  $1/2$-Lipschitz-continuous, then so is 
$f_1 \wedge f_2$. 
\item [(b)] If $f_1$ and $f_2$ are  piecewise continuously differentiable, 
then so is $f_1 \wedge f_2$. 
\end{enumerate}
\end{lem}
For the proof of these easy facts we refer to \cite{Ri2}. 
A function $f: \R^2 \to \R$ is called $1/2$-$\einh$-Lipschitz-continuous 
or piecewise continuously $\einh$-differentiable if for all  $r_2 \in \R$ 
the function  $f(.,r_2)$ is  $1/2$-Lipschitz-continuous or piecewise 
continuously differentiable respectively.
\begin{lem} \label{lemeig}
For $x' \in \R^2$ and $t \in \R$ the function 
$\tau_n(|.|) \wedge \mxt$ is  $1/2$-$\einh$-Lipschitz-continuous and 
piecewise continuously $\einh$-differentiable.
\end{lem}   
For details of the proof we again refer to \cite{Ri2}. Basically 
Lemma \ref{lemeig} follows from Lemma \ref{leanalmin}. The only 
difficulty is to show the continuity of $\tau_n(|.|) \wedge \mxt$, 
which might be a problem because of  the jump to infinity of $\mxt$ 
in case of $\hxt \cf \le 1/2$. But if 
$x \in \partial \{\mxt < \infty\} = \partial  \{ \fc(. -x') < 1 \}$
then $x-x'$ is contained in the closure of $\Kep$. Hence  $|x - x'| \le \cK$, 
which implies $|x'|-\cK \le |x|$. As $\tn$ is decreasing we obtain 
\[
\tau_n(|x|) \, \le \, \tn (|x'|-\cK) \, 
\le \, t + \hxt \, \le \, \mxt(x) 
\]
by definition of $\hxt$, which implies the claimed continuity. 

%=====================================================================

\subsection{Properties of the construction: Lemma \ref{leconshc}} 
\label{secleconshc}

We will first investigate monotonicity and regularity properties 
of $\etn^k$ and $\eTn^k$:     
\begin{lem} \label{leliphc}
For $X \in \X$ and $k \ge 0$
\begin{align}
&\etn^k \text{ is 1/2-$\einh$-Lipschitz-continuous and piecewise 
cont. $\einh$-differentiable}, \label{liphc} \\ 
&\eTn^k \text{ is $\lear$-increasing and bijective}. \label{wabihc}
\end{align} 
\end{lem}
\Bew
$\etn^k$ is the minimum of finitely many functions of the form 
$\tau_n(|.|) \wedge \mxt$, where $x' \in \R^2$ and $t \in \R$. 
Hence \eqref{liphc} is an immediate consequence of Lemmas 
\ref{lemeig} and \ref{leanalmin}. \eqref{liphc} implies that $\eTn^k$ is 
$\einh$-continuous and $\lear$-increasing, and hence bijective. 
This shows \eqref{wabihc}. \qed\\

\noindent
For \eqref{monohc} it suffices to observe that for every  $2 \le k \le m$ 
we have 
\[
\etaun^k \, = \, \etn^k(\epn^k) \, 
= \, \etn^{k-1}(\epn^k) \wedge m_{\epn^{k-1},\etaun^{k-1}}(\epn^k) \, 
\ge \, \etaun^{k-1}.
\]
This follows from the definition of $\etaun^k$ and $\etn^k$, 
from $\etn^{k-1}(\epn^k) \ge \etaun^{k-1}$ by the definition
of $\epn^{k-1}$ and from $\mxt \ge t$.\\ 
For \eqref{verconhc} and \eqref{grepshc} let $x_1,x_2 \in X$. 
Without loss of generality we may suppose that 
$x_1=\epn^j$ and $x_2=\epn^i$, where  $0 \le i \le j$. Here $\epn^0$ is 
interpreted to be any point of $X_{\Lan^c}$. We first observe that 
$\epn^j \in \ela^i := \{x \in \R^2: \etn^i(x) \ge \etaun^i\}$ and  
\begin{equation} \label{Klahc}
\alle x \in (\epn^i + \Ko) \cap \ela^i: \quad  
\etn^{j}(x) \, =\, \etn^{i}(x) \wedge \bigwedge_{i \le k \le  j} \, 
 m_{\epn^k,\etaun^k}(x) \, = \, \etaun^i.
\end{equation}
This holds as $\etn^i(x) \ge \etaun^i$ by definition of $\ela^i$, 
$m_{\epn^k,\etaun^k} \ge \etaun^i$ by \eqref{monohc} and
$m_{\epn^i,\etaun^i}(x) = \etaun^i$ by $x \in \epn^i + \Ko$.
If $\epn^j-\epn^i \in \Ko$, then $\epn^j \in  (\epn^i + \Ko) \cap \ela^i$, 
so \eqref{Klahc} implies   
$\etaun^j = \etn^{j}(\epn^j) = \etaun^i$, which shows \eqref{verconhc}.
For \eqref{grepshc} suppose $\epn^j-\epn^i \notin \Ko$. We have 
$\epn^j \in  \ela^i \weg (\epn^i + \Ko)$ and 
$\etaun^j = \etn^{j}(\epn^j)$ by definition, so it suffices to show
\begin{equation}  \label{grepsthc} 
\eTn^{j}(\ela^i \weg (\epn^i + \Ko)) \,
= \, \ela^i  \weg (\epn^i + \Ko)  +  \etaun^i \einh.  
\end{equation}
In order to show this we fix $r \in \R$. Continuity of $\etn^{i}(.,r)$ 
implies $\etn^i = \etaun^i$ on $\partial \ela^i(.,r)$. Just as in the proof 
of \eqref{Klahc} it follows that $\etn^{j} = \etaun^i$ on $\partial \ela^i(.,r)$. 
But $\eTn^{j}(.,r)$ is increasing, continuous and bijective by \eqref{wabihc}, 
so 
\[
\eTn^{j}(\ela^i) \, = \, \ela^i + \etaun^i \einh,
\] 
and combining this with  \eqref{Klahc} we are done.

%============================================================================

\subsection{Properties of the deformed translation: Lemma~\ref{leinnenaussenhc}}
\label{secpdt}

The following lemma shows how to estimate the translation distances $\etaun^k$.
\begin{lem} \label{leeigkonshc}
For $X \in \X$ and $k \ge 0$ we have 
\begin{align}
&\etaun^k \, \le \, \etn^0(\epn^{k'}) \quad 
\text{ for all }  k' \ge k 
, \label{tauobenhc}\\
&\etaun^k \, \ge \, \etn^0(\anx(\epn^k))    \quad \text{ if } 
  X \in \Gn.  \label{tauuntenhc}
\end{align} 
\end{lem}
\Bew
\eqref{tauobenhc} follows from the definition of $\epn^k$ and from 
$\etn^k \le \etn^0$.
For the proof of \eqref{tauuntenhc} let $X \in \Gn$. We first would like 
to show that 
\begin{equation} \label{verzerhc}
\forall x, x' \in X: \;
 |x| \le |x'|, \, x \stackrel{X,\Kepn}{\longleftrightarrow} x' \; 
 \Rightarrow  \; |\tn(|x|-\cK) - \tn(|x'|)| \cf \le 1/2. 
\end{equation}
Defining
\begin{equation} \label{sig1hc} 
\Sigma_1(n,X)\, :=  \sum_{x,x' \in X} 1_{\{|x| \le |x'|\}} 
 1_{\{ x \stackrel{X,\Kepn}{\longleftrightarrow} x'\}} 4
 \big(\tn(|x|-\cK) - \tn(|x'|)\big)^2 \cf^2
\end{equation}
we have  $\Sigma_1(n,X) <  1$ by definition of the set $\Gn$ of good 
configurations in \eqref{goodhc} and by  $X \in \Gn$. Hence every summand of 
$\Sigma_1$ is less than $1$, which implies \eqref{verzerhc}. We now can prove
\eqref{tauuntenhc} by induction on $k$.  For  $k = 0$ we have equality if 
the right hand side is defined to be $0$. For the inductive 
step $k-1 \to k$ let  $i \le k-1$.  By  \eqref{tauobenhc}, 
the inductive hypothesis and  \eqref{verzerhc} we have 
\[
0 \, \le \, \big(\tn(|\epn^i|-\cK) - \etaun^{i}\big) \cf  \,  \le \, 
\big(\tn(|\epn^i| - \cK) - \tn(|\anx(\epn^i)|)\big) \cf \, 
\le \, 1 / 2, 
\]
so $h_{\epn^i,\etaun^i} \cf \le 1/2$. 
Therefore $m_{\epn^i,\etaun^i}(\epn^{k}) = \infty$ whenever 
$\epn^k-\epn^i \notin \Kep$. Thus 
\[ \etaun^{k} \, =\, \etn^{k}(\epn^{k})\, 
 =  \, \etn^0(\epn^{k}) \wedge \bigwedge_{i < k: \epn^{k}-\epn^i \in \Kep } 
 m_{\epn^i,\etaun^i}(\epn^{k}) \, \ge \, \etn^0(\anx(\epn^k)),
\]
where the last step follows from 
$m_{\epn^i,\etaun^i}(\epn^{k}) \ge \etaun^i \ge \etn^0(\anx(\epn^i))$, 
which holds by  induction hypothesis, and 
$\anx(\epn^i)=\anx(\epn^k)$ for $\epn^{k}-\epn^i \in \Kep$. \qed\\

\noindent
In the proof of \eqref{tauuntenhc} we have also shown that 
good configurations $X \in \Gn$ have the following property: 
In the construction of $\Tn(X)$ we have 
$h_{\epn^k,\etaun^k}\cf \le 1/2$ for every $k$, i.e. 
in the definition of $m_{\epn^k,\etaun^k}$ we always have the second case. 
Now we will prove Lemma \ref{leinnenaussenhc}. It suffices to show for all 
$X \in \Gn$ and $x \in X$ that  
\begin{equation} \label{inaushc}
\begin{split}
&x \in \La_{n'} \; \Rightarrow  \; \tnx(x) = \tau, \qquad 
x \in \La_{n'}^c \; \Rightarrow \; x + \tnx(x) \einh - \tau \einh \notin \La_{n'-1}\\
&x \in \Lan^c  \; \Rightarrow \;  \tnx(x) = 0, \qquad \text{ and } \quad  
x \in \Lan  \; \Rightarrow \; x + \tnx(x)\einh \in \Lan. 
\end{split}
\end{equation}
So let $X \in \Gn$ and $x \in X$. We first note that 
\begin{equation}  \label{abdeftxhc}
0 \, \le \, \tn(|\anx(x)|) \, 
\le \, \tnx(x) \, \le \, \tn(|x|) \, \le \, \tau,
\end{equation}
which is an immediate consequence of \eqref{tauuntenhc} and \eqref{tauobenhc}.
We observe 
\[
x \in \La_{n'} \; \Rightarrow \; \anx(x) \in \La_R \; \Rightarrow  \; 
\tn(|\anx(x)|) = \tau  \; \Rightarrow \; \tnx(x) = \tau,
\]
where we have used the definition of $R$, $X \in \Gna$, \eqref{tauinnenaussenhc} 
and \eqref{abdeftxhc}. This gives the first assertion of \eqref{inaushc}. 
The second assertion is an immediate consequence of $0 \le \tau- \tnx(x) \le 1$, 
which follows from \eqref{abdeftxhc} and $\tau \le 1$. 
The third assertion follows from \eqref{abdeftxhc} and \eqref{tauinnenaussenhc}, 
and for the fourth assertion let  $x \in \Lan$. As 
\[
x  \lear  x + \tnx(x)\einh  \lear  \eTn^0(x)
\]
by  \eqref{abdeftxhc}, it suffices to show that also $\eTn^0(x) \in \Lan$. 
This however follows from   $\eTn^0 = id$ on ${\Lan}^c$ and the bijectivity 
of  $\eTn^0$ from \eqref{wabihc}. 
%==========================================================================

\subsection{Bijectivity of  the transformation: Lemma \ref{lebijhc}}
\label{secbijhc}

We will construct the inverse transformation $\tTn$ recursively 
just as in the construction of $\Tn$, 
i.e. from a given configuration $\tX$ we  will choose points $\etpn^k$ 
and translate them by $\ettaun^k$ in direction $-\einh$.
\begin{figure}[!htb] 
\begin{center}
\psfrag{1}{$\epn^0$}
\psfrag{2}{$\etpn^0$}
\psfrag{5}{$\epn^2$}
\psfrag{6}{$\etpn^2$}
\psfrag{7}{$\epn^1$}
\psfrag{8}{$\etpn^1$}
\psfrag{A}{$\etaun^0$}
\psfrag{B}{$\ettaun^0$}
\psfrag{E}{$\etaun^2$}
\psfrag{F}{$\ettaun^2$}
\psfrag{G}{$\etaun^1$}
\psfrag{H}{$\ettaun^1$}
\psfrag{x}{$\Tn: X \mapsto \tX$}
\psfrag{y}{$\tTn: \tX \mapsto X$}
\includegraphics[scale=0.5]{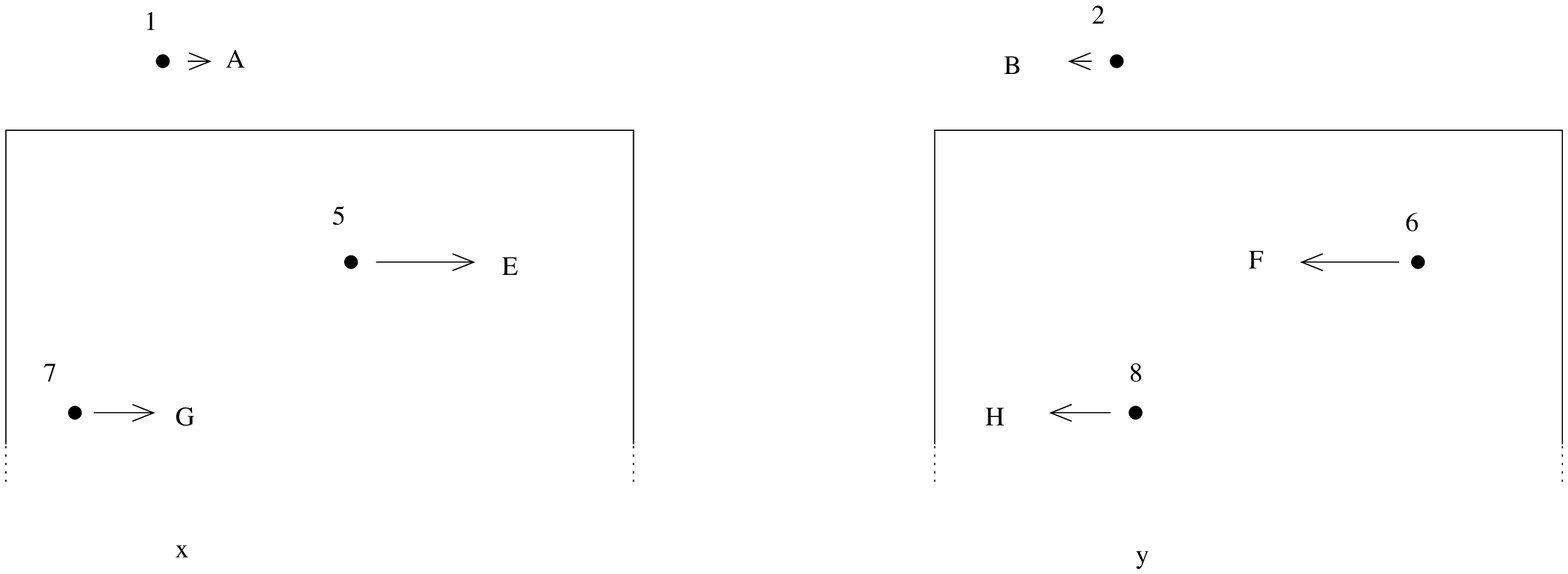}
\end{center}
\caption{Construction of the inverse  $\tTn$ of $\Tn$.} \label{figinvhc}
\end{figure}

\noindent
To get an idea how to define the inverse transformation we start with 
$X \in \X$ and set  $\tX:= \Tn(X)$. In the construction of $\tX$ we defined  
points $\epn^k$ and translation distances $\etaun^k$. We denote the 
corresponding image points by $\etpn^k := \epn^k + \etaun^k \einh$, 
see Figure~\ref{figinvhc}. For the construction of the inverse transformation 
we have to find a method to identify the points $\etpn^k$ among the points 
of $\tX$ without knowing $X$. Suppose now that we have already found 
$\etpn^1, \ldots , \etpn^{k-1}$. Then inductively we are able to construct
the translation distances $\etaun^{i}$ for all $1 \le i < k$, because
$\etn^i$ is defined in terms of $\epn^j$ and $\etaun^j$ where $j<i$,
$\eTn^i = id + \etn^i \einh$, $\epn^i = (\eTn^i)^{-1}(\etpn^i)$ and 
$\etaun^i = \etn^i(\epn^i)$. So in particular we know the 
transformation functions $\etn^k$ and $\eTn^k$.  
Thus the following lemma gives a characterisation of $\etpn^k$
just as needed: 
\begin{lem} \label{leminimhc}
Let $1 \le k \le m$. For every   
$\tx \in \tX_{\Lan} \weg \{\etpn^1, \ldots, \etpn^{k-1}\}$ we have  
\[
\etn^k \circ (\eTn^k)^{-1}(\etpn^k) \,\le \,\etn^k \circ (\eTn^k)^{-1}(\tx).
\] 
For all $\tx$ for which equality occurs we have  
$(\eTn^k)^{-1} (\etpn^k) \le  (\eTn^k)^{-1}(\tx)$. 
\end{lem} 
\Bew 
We first observe that for all $k$  by definition of $\eTn^k$ we have
\begin{equation} \label{trdtinvhc}
(\eTn^k)^{-1} + \etn^k \circ (\eTn^k)^{-1}\einh \, = \, id.  
\end{equation}
Since $\etn^{k+1} \le \etn^{k}$, we also have $\eTn^{k+1} \lear \eTn^{k}$, 
and therefore  $(\eTn^{k})^{-1} \lear  (\eTn^{k+1})^{-1}$ by the 
$\einh$-monotonicity  of $(\eTn^{k+1})^{-1}$ from \eqref{wabihc}. 
Together with \eqref{trdtinvhc} this implies 
\begin{equation} \label{inveincrhc}
\etn^{k+1} \circ (\eTn^{k+1})^{-1} \, \le \, \etn^k \circ (\eTn^k)^{-1}.
\end{equation}
Now let $1 \le k \le m$ and   
$\tx \in \tX_{\Lan} \weg  \{\etpn^1, \ldots, \etpn^{k-1}\}$, i.e.  
$\tx = \etpn^l$ for some $l \ge k$. By definition we have 
$\etn^{l}(\epn^l)  =  \etaun^l$, $\eTn^{l}(\epn^l) = \etpn^l$ and 
$\etpn^k = \eTn^{k}(\epn^k)$. Using  
\eqref{monohc} and \eqref{inveincrhc} we deduce 
\[
 \etn^{k} (\eTn^{k})^{-1} (\etpn^k)\, 
= \, \etaun^k\, \le \, \etaun^l \, 
= \, \etn^{l} (\epn^l) \, 
= \, \etn^{l} (\eTn^{l})^{-1} (\tx) \,
\le \, \etn^{k} (\eTn^{k})^{-1} (\tx).
\]
If for the given $\tx$ we have equality, all inequalities in the previous line 
have to be equalities, so $\etaun^k = \etaun^l$  and 
$\etn^{l} (\eTn^{l})^{-1} (\tx) =\etn^{k} (\eTn^{k})^{-1} (\tx)$. 
Combining this with \eqref{trdtinvhc} we get 
$\epn^l = (\eTn^{l})^{-1} (\tx) =(\eTn^{k})^{-1} (\tx)$, so 
$\eTn^k(\epn^l) = \tx$, and thus $\etn^k(\epn^l) = \etaun^l = \etaun^k$. 
By definition of $\epn^k$ we conclude 
$(\eTn^k)^{-1} (\etpn^k) = \epn^k \le \epn^l$ and we are done.  \qed \\

\noindent
Lemma \ref{leminimhc} tells us exactly how to construct the inverse of $\Tn$ 
recursively. So let $\tX \in \X$.
Let $\tilde{m}=\tilde{m}(\tX) := \# \tX_{\Lan}$,  
$\ttnx^0 = \tau_n(|.|)$  and $\ttaunx^0 := 0$. 
In the k-th construction step $(1 \le k \le \tilde{m})$ let
\[
\ttnx^{k} \, :=  \, \ttnx^{k-1} \wedge 
   m_{\tpnx^{k-1}-\ttaunx^{k-1},\ttaunx^{k-1}}, \quad \text{ where } \; 
m_{\tpnx^{0}-\ttaunx^{0},\ttaunx^{0}} \,
:= \, \bigwedge_{\tx \in \tX_{\Lan^c}} m_{\tx,0}.
\]
Let $\tTnxn^k = id + \ttnx^k \einh$ and let $\tpnx^k$ be the point of 
$\tX_{\Lan} \weg \{\tpnx^1, \ldots, \tpnx^{k-1}\}$ at which the minimum
of $\ttnx^k \circ (\tTnxn^k)^{-1}$ is attained. If there is more than one 
such point then take the point $y$ such that $(\tTnxn^k)^{-1}(y)$ is minimal 
with respect to the lexicographic order $\le$. Let  
$\ttaunx^k := \ttnx^k \circ (\tTnxn^k)^{-1}(\tpnx^k)$ be  the corresponding 
minimal value.
In the above notations we will omit dependencies on $\tX$
if it is clear which configuration is considered. 
We need to show that the above construction is well defined, i.e. 
that   $\tTnxn^k$ is invertible in every step. Furthermore we 
need some more properties of the construction: 
\begin{lem} \label{leeiginvhc}
Let $\tX \in \X$ and $k \ge 0$. Then
\begin{align}
& \ettn^k \text{ is 1/2-$\einh$-Lipschitz-continuous}, \; 
\etTn^k \text{ is bijective and $\lear$-increasing}, \label{lipinvhc}\\
& (\etTn^k)^{-1} + \ettn^k \circ (\etTn^k)^{-1} \einh  \, = \,id, 
 \label{TSinvhc}\\
& \alle c \in \R ,x \in \R^2: \; \ettn^k \circ (\etTn^k)^{-1}(x) 
 \,\ge\,c\, 
  \Leftrightarrow \, \ettn^k(x -  c\einh) \, \ge \, c, \label{tsinvhc}\\ 
& \ettn^{k} \, \le \, \ettn^{k-1} \quad \text{ and } \quad  
 \ettaun^{k-1} \, \le \, \ettaun^{k}. \label{monoinvhc}
\end{align} 
\end{lem}
\Bew
The definitions of $\ettn^{k}$ and $\etTn^k$ are similar to those of 
$\etn^{k}$ and $\eTn^k$, so we can show \eqref{lipinvhc} and \eqref{TSinvhc} 
just as the corresponding properties in \eqref{liphc}, \eqref{wabihc}  and  
\eqref{trdtinvhc}. 
For \eqref{tsinvhc} we note that for $c \in \R$ and $x \in \R^2$ 
the equivalence 
\[
\begin{split}
\ettn^k \circ (\etTn^k)^{-1}(x) \,\ge \,c \quad &\Leftrightarrow \quad  
(\etTn^k)^{-1}(x)  \, \lear \, x - c\einh \\ 
&\Leftrightarrow \quad  
x \, \lear \, \etTn^k(x - c\einh ) = x - c\einh + \ettn^k(x - c\einh )\einh
\end{split}
\]
follows from   \eqref{TSinvhc} and \eqref{lipinvhc}.
The first part of \eqref{monoinvhc} is obvious and for the second part  
we observe that 
\[
\begin{split}
\ettn^{k-1} &\circ (\etTn^{k-1})^{-1}(\etpn^k) \ge \ettaun^{k-1} 
\quad \Rightarrow \quad 
\ettn^{k-1}(\etpn^k - \ettaun^{k-1}\einh ) \ge \ettaun^{k-1}  \\ 
&\Rightarrow \quad 
\ettn^{k}(\etpn^k - \ettaun^{k-1} \einh) \ge \ettaun^{k-1}  
\quad \Rightarrow \quad 
\ettaun^k =\ettn^{k} \circ (\etTn^k)^{-1}(\etpn^k) \ge \ettaun^{k-1},
\end{split} 
\]
where the first statement holds by definition of $\etpn^{k-1}$, 
the first and the third implication hold by  \eqref{tsinvhc} and the 
second holds by definition of $\ettn^k$.
\qed\\ 

\noindent
Let $\ttnx(\tpnx^k) := \ttaunx^{k}$ and $\ttnx(x)=0$ for $x \in \tX_{\Lan^c}$. 
This defines a translation distance function $\ttnx: \tX \to \R$. 
Let $\tTn: \X  \to \X$ be defined by 
\[ 
\tTn(\tX) \, 
:= \, \tX_{\Lan^c} \cup \{\tpnx^k - \ttaunx^k \einh: 1 \le k \le m\} 
= \{ x - \ttnx(x) \einh : x \in \tX \}.
\]
By Lemma \ref{lemeas} we again see that all above objects are 
measurable with respect to the considered $\si$-algebras. 
The only difficulty is to show that the functions   $(\tTnxn^k)^{-1}(x)$  
are measurable, which follows from the  $\einh$-monotonicity of  
$\tTnxn^k$.\\

In order to show that  $\tTn$ really is the inverse of $\Tn$ we need an 
analogue of Lemma \ref{leminimhc}. Let  $\tX \in \X$. 
Let $\ettn^k$, $\etTn^k$, $\etpn^k$ and  $\ettaun^k$ 
$~(0 \le k \le \tilde{m})$ as above and denote $X := \tTn(\tX)$ and 
$\epn^k := \etpn^k - \ettaun^k \einh$, see Figure \ref{figinvhc}.
\begin{lem} \label{leminiminvhc}
Let $1 \le k \le \tilde{m}$. For every 
$x \in X_{\Lan} \weg \{\epn^1, \ldots, \epn^{k-1}\}$ we have 
\[
\ettn^k (\epn^k) \le \ettn^k (x).
\] 
For all $x$ for which equality occurs we have $\epn^k \le x$.  
\end{lem} 
\Bew 
Let $1 \le k \le \tilde{m}$ and 
$x \in X_{\Lan} \weg \{\epn^1, \ldots, \epn^{k-1}\}$, i.e.  
$x = \epn^l$ for some $l \ge k$. By definition of $\ettaun^k$ and 
$\ettaun^l$ and \eqref{TSinvhc} we have  $(\etTn^{k})^{-1}(\etpn^k) = \epn^k$
and $(\etTn^{l})^{-1}(\etpn^l) = x$. Using \eqref{monoinvhc} we obtain  
\[
\ettn^{k} (\epn^k)\, 
= \,  \ettaun^k\, \le \, \ettaun^l \, 
= \, \ettn^{l} (\etTn^{l})^{-1} (\etpn^l) \,
= \, \ettn^{l} (x) \, 
\le \, \ettn^{k} (x).
\]
If for the given $x$ we have equality, all inequalities in the previous line
have to be equalities, so $\ettaun^k = \ettaun^l$ and 
$\ettn^{k} (x) =\ettaun^l$, i.e.   
$\etTn^{k} (x) = x +  \ettaun^l \einh = \etpn^l$. 
This gives $\ettaun^k = \ettaun^l =  \ettn^{k} (x) 
=  \ettn^{k} \circ (\etTn^{k})^{-1} (\etpn^l)$. So 
$\epn^k =  (\etTn^{k})^{-1}(\etpn^k) \le  (\etTn^{k})^{-1} (\etpn^l) = x$ 
by definition of $\etpn^k$ and we are done. \qed 
\begin{lem}\label{leinvhc}
On $\X$ we have $ \quad  \tTn \circ \Tn \, = \, id \quad $ 
and  $ \quad  \Tn \circ \tTn \, = \, id$.
\end{lem}
\Bew
For the first part let  $X \in \X$ and $\tX := \Tn(X)$. We have 
$\tilde{m}(\tX) = m(X)$ by construction and we have $X_{\Lan^c} = \tX_{\Lan^c}$
by \eqref{innenaussenhc}. Now 
it suffices to prove 
\begin{equation}  \label{tildegleichhc}
\ttnx^k \, = \, \tnx^k, \; \tTnxn^k \, = \, \Tnxn^k, \; 
\ttaunx^k \, = \, \taunx^k \, \text{ and } \, 
\tpnx^k \, = \, \pnx^k +  \taunx^k  
\end{equation} 
for every $k \ge 0$ by induction on  $k$. Here $\tpnx^0 = \pnx^0 +  \taunx^0$
is interpreted as $X_{\Lan^c} = \tX_{\Lan^c}$. The case $k=0$ is trivial. 
For the inductive step  $k-1 \to k$ we observe that  
$\ettn^k = \etn^k$ by induction hypothesis, and $\etTn^k  =  \eTn^k$ 
is an immediate consequence. Combining this with Lemma \ref{leminimhc}
and the definition of $\etpn^k$ we get 
$\etpn^k = \epn^k  +  \etaun^k$ and   $\ettaun^k = \etaun^k$.\\
For the second part let  $\tX \in \X$ and 
$X := \tTn(\tX)$. As above it suffices to show \eqref{tildegleichhc}
by induction on $k$. Here $\tX_{\Lan^c} = X_{\Lan}$ follows from an analogue 
of \eqref{innenaussenhc} and the inductive step follows from Lemma  
\ref{leminiminvhc}.
\qed
  
%===========================================================================

\subsection{Density of the transformed process: Lemma \ref{ledensphc}}
\label{secledensphc}

By definition the left hand side of \eqref{denshc} equals 
\[ 
e^{-4 n^2} \sum_{k \ge 0} \frac{1}{k!} I(k), \quad
\text{ where }  \,
  I(k) \, = \, \int_{{\Lan}^{k}} dx (f \circ \Tn \cdot \ph) (\bX_x),
\]
using the shorthand notation $\bX_x := \{ x_i: i \in J\} \cup \bX_{\Lan^c}$
for  $x  \in \Lan^{\;J}$. To compute $I(k)$ we need to calculate  
$\Tn(\bX_x)$, and for this we must identify the points $P^i_{n,\bX_x}$ among 
the particles $x_j$. So let $\Pi$ be the set of all permutations 
$\per:\{1,\ldots ,k\}\to \{1,\ldots,k\}$. For $\per \in \Pi$ 
let 
\[ \begin{split}
\ak \, &:= \, \big\{ x \in\Lan^{\;k}: \alle 1 \le j \le k: 
 x_{\per(j)} = P^j_{n,\bX_x}\big\}\quad \text{ and } \\
\tak \, &:= \, \big\{ x \in\Lan^{\;k}: \alle 1 \le j \le k:
 x_{\per(j)} = \tilde{P}^j_{n,\bX_x}\big\},
\end{split}\]
where  $\tilde{P}^j_{n,\bX_x}$ are the points from the construction 
of the inverse transformation in Subsection \ref{secbijhc}.  
Now we can write
\[
I(k) \, = \,\sum_{\per \in \Pi} I(k,\per), \quad \text{ where } \,
 I(k,\per)\, = \, \int_{{\Lan}^{k}} dx \, 
        1_{\ak}(x) (f \circ \Tn \cdot \ph)(\bX_x).   
\]
If  $x \in \ak$ we can derive a simple expression for $\Tn(\bX_x)$:
For $x \in \Lan^{\;k}$ we define a formal transformation 
$T_{\per}(x)  :=  (T^i_{\per,x} (x_i))_{1 \le i \le k}$, where 
\[
T^{\per(j)}_{\per,x} :=  id + t^{\per(j)}_{\per,x}\einh, \quad 
t^{\per(j)}_{\per,x}  :=  t_{n,\bX_{x^{\per,j-1}}}^{j} \quad \text{ and } \quad  
x^{\per,j-1}  :=  (x_{\per(i)})_{1 \le i \le j-1}.
\]
Clearly, $T^{\per(j)}_{\per,x}$ doesn't depend
on all components of $x$, but only on those $x_{\per(l)}$ such that 
$l \le j-1$. By definition we now have 
\begin{equation} \label{athc}
\alle x \in \ak: \quad \T_n(\bX_x) \, = \, \bX_{T_{\per}(x)} \; \text{ and } \;
\, T^j_{n,\bX_x} \, = \, T^{\per(j)}_{\per,x} \, \text{ for all } j \le k.
\end{equation}
Furthermore we observe that for all  $x \in (\R^2)^{k}$  we have 
\begin{equation} \label{trapohc}  
x \in \ak \quad \Leftrightarrow \quad T_{\per}(x) \in \tak.
\end{equation}
Here ``$\Rightarrow $'' holds by \eqref{athc} and \eqref{tildegleichhc} 
from the proof of Lemma  \ref{leinvhc}. For  ``$\Leftarrow $''
let  $x \in  (\R^2)^{k}$ such that $T_{\per}(x) \in \tak$ and
let $X':= \tTn(\bX_{T_{\per}(x)})$, 
where $\tTn$  is the inverse of  $\Tn$ as defined in the last subsection. 
By induction on $j$ we can show 
\[
\alle 1 \le j \le k:  \quad T_{n,X'}^j =  T^{\per(j)}_{\per,x} \quad \text{ and } \quad  
x_{\per(j)} = P_{n,X'}^j.
\]
In the inductive step $j-1 \to j$ the first assertion follows from the
induction hypothesis and the second 
follows from the bijectivity of $T_{n,X'}^j$ and  
\[
T_{n,X'}^j(x_{\per(j)}) \, = \,  T^{\per(j)}_{\per,x}(x_{\per(j)}) \, 
= \, \tilde{P}^j_{n,\bX_{T_{\per}(x)}} \, 
= \, P^j_{n,X'} + \tau^j_{n,X'} \, = \, T_{n,X'}^j(P_{n,X'}^j), 
\]
which follows from $T_{n,X'}^j =  T^{\per(j)}_{\per,x}$, the definition of $\tak$ 
and \eqref{tildegleichhc} from the proof of Lemma \ref{leinvhc}. 
This completes the proof of the above assertion 
and we conclude $\bX_x= X'$, which implies 
$x_{\per(j)}  = P_{n,X'}^j = P^j_{n,\bX_x}$. Thus \eqref{trapohc} holds. \\
Now let $g: (\R^2)^{k} \to \R$, $g(x) :=  1_{\tak}(x) f(\bX_x)$.
Then \eqref{athc} and \eqref{trapohc} imply 
\[
I(k,\per) \, = \,  \Big[ \prod_{j=1}^{k}  \int dx_{\per(j)} \,
 \big|1 + \partial_{1} t^{\per(j)}_{\per,x}(x_{\per(j)})\big| \Big] \, 
 g(T_{\per}(x)),  
\]
where we have also inserted the definition  of $\ph$ \eqref{densdefhc}.
Now we transform the integrals. For $j=k$ to $1$ we substitute 
$x_{i}' := T^{i}_{\per,x} x_{i}$, where $i := \per(j)$. 
The transformation only concerns the first component of 
$x_i = (r_i,\bar{r}_i)$. For fixed $\bar{r}_{i}$ $~r_{i}$ is transformed by 
$id +  t^{i}_{\per,x}(.,\bar{r}_{i})$. From \eqref{liphc} we know 
that $t^{i}_{\per,x}(.,\bar{r}_{i})$ is $1/2$-Lipschitz-continuous 
and piecewise continuously differentiable, 
so  $id +  t^{i}_{\per,x}(.,\bar{r}_{i})$ 
is strictly increasing and piecewise continuously differentiable. 
Therefore the Lebesgue transformation theorem  \eqref{lebtrsatz} gives 
\[
dx_{i}' \, = \,  dx_{i} \big| 1 + \partial_{1} t^{i}_{\per,x}  
(x_{i})\big|.
\] 
Thus 
\[
I(k,\per) \, = \,  \Big[\prod_{j=1}^{k}  \int dx'_{\per(j)} \Big] \, g(x') \, 
=  \, \int_{{\Lan }^k} dx \, 1_{\tak}(x) f(\bX_x),
\]
and we are done as the same arguments show that the right hand side of 
\eqref{denshc} equals 
\[ 
 e^{-4 n^2} \sum_{k \ge 0} \frac{1}{k!} \sum_{\per \in \Pi} 
 \int_{{\Lan }^k} dx \,   1_{\tak}(x) f(\bX_x).
\]
An analogous argument shows that the density function is well defined: 
\[
\begin{split} 
\alle \bX \in \X: \quad 
\nu_{\Lan} &(``\ph \text{ is well defined}''|\bX)\\ 
&= \,  e^{-4 n^2} \sum_{k \ge 0} \frac{1}{k!}
 \sum_{\per \in \Pi}   \int_{{\Lan }^k} dx \,  1_{\ak}(x)  \prod_{j=0}^{k} 
 1_{\{ \partial_{1} t^{\per(j)}_{\per,x}(x_{\per(j)}) \text{ exists} \}}.
\end{split}
\]
As $t^{\per(j)}_{\per,x}$ is piecewise continuously  $\einh$-differentiable,
we have for arbitrary  $r \in \R$, $k$, $\per$ and $x$ 
as above that $\partial_{1} t^{\per(j)}_{\per,x}(.,r)$ exists $\la^1$-a.s..
So we may replace all indicator functions in the above product by $1$ using
Fubini's theorem. Hence the above probability equals $1$. 

%=======================================================================

\subsection{Estimation of the densities: Lemma \ref{ngrosshc}}
\label{secngrosshc}

Let $X \in \Gn$. By the  $1/2$-$\einh$-Lipschitz-continuity from 
\eqref{liphc} we have 
\[
|\partial_{1} \tnx^k(\pnx^k)| \le 1/2.
\]
Using $-\log(1-a) \le 2a$ for $0 \le a \le 1/2$ we obtain
\[ \begin{split} 
 f_n(X) \, &:= - \log \bph(X)  -  \log \ph(X) \\
&= \, - \sum_{1 \le k \le m} 
 \log\big(1 - (\partial_{1} \tnx^k(\pnx^k))^2\big) \, 
 \le \, \sum_{1 \le k \le m} 2 (\partial_{1} \tnx^k(\pnx^k))^2.
\end{split} \]
If $\partial_{1} \etn^k(\epn^k)$ exists it equals either 
$\partial_{1} \etn^0(\epn^k)$ or  $\partial_1 m_{x,\etn(x)}(\epn^k)$ for some  
$x \in X$ such that $x \neq \epn^k$ and $\epn^k \in x + \Kep$. 
By using \eqref{tauuntenhc} we see that 
\begin{displaymath}
|\partial_1 m_{x,\etn(x)}(\epn^k)| \, 
\le \, (\tn(|x|-\cK) -\etn(x)) \cf\,
\le \, \big(\tn(|x|-\cK) -\tn(|\anx(x)|)\big) \cf. 
\end{displaymath}
Furthermore $|\partial_{1} \etn^0(\epn^k)| \le \tau q(|\epn^k| - R)/Q(n-R)$
by definition of $\etn^0 = \tn(|.|)$,
so we can estimate $f_n(X)$ by the sum of the two following terms: 
\begin{equation} \begin{split} \label{sig2hc} 
\Sigma_2(n,X) \,  &:= \,  2\tau^2 \sum_{x \in X}
  1_{\{x \in \Lan\}} \frac{q(|x| -R)^2}{Q(n-R)^2},  \\
\Sigma_3(n,X) \, &:=   \, 2 \cf^2 \sideset{}{^{\neq}} \sum_{x,x' \in X} 
 \sideset{}{} \sum \limits_{x'' \in X} 1_{\Kep}(x'-x)
   1_{\{x \stackrel{X,\Kepn}{\longleftrightarrow} x''\}}  1_{\{|x|\le |x''|\}}\\
&\qquad \hspace{3 cm} \times (\tn(|x|-\cK) -\tn(|x''|))^2.
\end{split} \end{equation}
Using these terms in the definition \eqref{goodhc} of $\Gn$ we are done.

%=====================================================================

\subsection{Set of good configurations: Lemma \ref{lebadsmallhc}}
\label{secsighc}

The functions $\Sigma_i(n,X)$ from the definition of the set of 
good configurations $\Gn$ in \eqref{goodhc} have been specified in 
\eqref{sig1hc} and \eqref{sig2hc}. Using the shorthand 
\[
\Taun(x,x'') \, := \, 1_{\{|x|\le |x''|\}} |\tn(|x|-\cK) - \tn(|x''|)|^2
\]
we have 
\[ \begin{split}
& \Sigma_1  =   4 \cf^2 \sideset{}{^{}}\sum \limits_{x,x'' \in X} 
  1_{\{x \stackrel{X,\Kepn}{\longleftrightarrow} x''\}}
  \Taun(x,x''),\quad  
 \Sigma_2  =  2\tau^2 \sideset{}{}\sum \limits_{x \in X}  1_{\{x \in \Lan\}}
 \frac{q(|x| -R)^2}{Q(n-R)^2} \quad \text{ and }\\ 
& \Sigma_3  =  2 \cf^2  \sideset{}{^{\neq}} \sum \limits_{x,x' \in X} 
 \sideset{}{} \sum \limits_{x'' \in X} 1_{\Kep}(x'-x)
   1_{\{x \stackrel{X,\Kepn}{\longleftrightarrow} x''\}} \Taun(x,x''). 
\end{split} \]
We will show that the expectation of every $\Sigma_i$ can 
be made arbitrarily small when $n$ is chosen big enough. But
first we will give some relations used later.\\ 
Let $n \ge R+1$. For $s' > s$ such that $s'  > R$ and $s < n$ we have  
\begin{displaymath}
0 \,\le \, r(s-R,n-R) - r(s'-R,n-R) \, 
= \,  \int_{R \vee s}^{s' \wedge n}
 \frac{q(t-R )}{Q(n-R)} dt \, 
\le \,  (s'-s) \, \frac{q(s - R)}{Q(n-R)}
\end{displaymath}
by the monotonicity of  $q$. Defining $\bn := n + \cK$ and $\bR := R + \cK$ 
we thus have 
\begin{equation} 
\label{taugegenqhc}
\Taun(x,x') \,
\le \, 1_{\{x \in \La_{\bn} \}} \tau^2 \, (|x'|-|x|+ \cK)^2 \, 
    \frac{q(|x| - \bR)^2}{Q(\bn-\bR)^2} \quad \text{ for } x,x' \in \R^2, 
\end{equation} 
using the substitution $s' := |x'|$ and $s := |x| - \cK$. (If  
$s' \le R$ or $s \ge n$ then $\Taun(x,x') = 0$.)
The following relations will give us control over the relevant terms 
of the right hand side of \eqref{taugegenqhc}. We first observe that  
\begin{equation}\label{intqhc}
\int_{\La_{\bn}} dx  \, q(|x| -  \bR)^2 \,
\le \, 16\bR^2 \,+ \, 32 Q(\bn-\bR) \quad \text{ for } \bn \ge 2 \bR.    
\end{equation}
Indeed, writing $s := |x|$ we obtain
\begin{displaymath}
\begin{split}
\int_{\La_{\bn}} &dx  \, q(|x| - \bR)^2 \quad 
\le \quad \int_0^{2\bR} ds \, 8s \, 
  + \, \int_{\bR}^{\bn-\bR} ds \, 8(s+\bR) q(s)^2 \\
&\le  \, 16\bR^2 \,+ \, 32 \int_0^{\bn-\bR} q(s)ds  \quad 
\le \quad 16\bR^2 \,+ \, 32 Q(\bn-\bR).
\end{split}
\end{displaymath}
In the first step we used $q \le 1$, and in the second step 
$\bR \le s$ and $sq(s) \le 2$. We observe 
$\lim \limits_{n \to \infty}Q(n) = \infty$, which is a consequence of 
$\log \log n \le Q(n)$ for $n >1$. Therefore by   
\eqref{intqhc}
\begin{equation} \label{ntoinfhc}
\lim_{n \to \infty} c(n)  \, = \, 0 \quad \text{ for } \quad 
c(n) \, := \, \int_{\La_{\bn}} dx  \, \frac{q(|x| -  \bR)^2}{Q(\bn-\bR)^2}.
\end{equation} 
Finally, for $x_0,\ldots ,x_m \in \R^2$ such that $x_i-x_{i-1} \in \Kep$ 
we have $|x_i-x_{i-1}| \le \cK$, so 
\begin{equation} \label{sqmaxhc}
(|x_m|-|x_0| + \cK)^2 \, 
    \le \,  (m+1)^2 \cK^2. 
\end{equation}
Now we will use the ideas of the proof of Lemma \ref{lereichhc}. 
For $X \in \X$ we can estimate the summands of 
$\Sigma_1(n,X)$ by considering all paths $x_0,\ldots ,x_m$ 
in the graph $(X,\Kepn)$ connecting $x=x_0$ and $x''=x_m$. By \eqref{taugegenqhc} 
and \eqref{sqmaxhc} we can estimate $\Sigma_1(n,X)$ by a constant $c$ times 
\[
\sum_{m \ge 0} (m+1)^2  \quad \sideset{}{^{\neq}}
 \sum_{x_0,\ldots ,x_m \in X}   1_{\{x_0 \in \La_{\bn} \}}
 \frac{q(|x_0|-\bR)^2}{Q(\bn-\bR)^2} 
\prod_{i=1}^m 1_{\{x_i x_{i-1} \in \Kepn\}}.
\]
Using Lemma~\ref{lekorab} we can thus proceed
as in the proof of Lemma~\ref{lereichhc}:
\[ 
\int \mu(dX) \, \Sigma_1(n,X)\,
\le  \, z\xi c \sum_{m \ge 0}  (m+1)^2 (\ex z \xi )^{m} c(n).
\] 
Likewise,
\begin{displaymath} 
\int \mu(dX) \Sigma_2(n,X)\,
\le \, 2 z \xi  \tau^2 c(n).
\end{displaymath}
Finally, we can estimate $\Sigma_3(n,X)$ by a constant $c$ times  
\[ 
\begin{split}
&\sum_{m \ge 0} (m+1)^2  \quad \sideset{}{^{\neq}}
 \sum_{x_0,\ldots ,x_m \in X}   1_{\{x_0 \in \La_{\bn} \}}
 \frac{q(|x_0|-\bR)^2}{Q(\bn-\bR)^2}\\
&\qquad \times \prod_{i=1}^m 1_{\{x_i x_{i-1} \in \Kepn\}} 
 \Big[ \sum_{x' \in X , x' \neq x_i \alle i} 1_{\Kep}(x'-x_0)
 +  \sum_{j=1}^{m}   1_{\Kep}(x_j-x_0) \Big]. 
\end{split} \]
The second sum in the brackets can be estimated by $m$. As above 
\[
\int \mu(dX)\, \Sigma_3(n,X) \, 
\le \, z \xi c \sum_{m \ge 0}  (m+1)^2 (\ex z \xi )^{m}  c(n)(z \xi \ex + m). 
\]
In the bounds on the expectations of $\Sigma_1$ and $\Sigma_3$ 
the sums over $m$ are finite   by \eqref{dechc}. 
Collecting all estimates and using \eqref{ntoinfhc} we thus find that   
\begin{displaymath}
\int \mu (dX) \sum_{i=1}^3 \Sigma_i(n,X) \, 
\le \, \frac{\de}{2}
\end{displaymath}
for sufficiently large $n$, and $\mu  (\Gn^c) \le \de$ follows from 
the high probability of  $\Gna$, the Chebyshev inequality and the definition of 
$\Gn$ in \eqref{goodhc}.

%============================================================================= 
%
%============================================================================= 

\section{Proof of Theorem \ref{sym}: Main steps} \label{secproofsym}

\subsection{Basic constants} \label{constants}

Let $(U,z,\XX)$ be admissible with Ruelle bound $\xi$, where 
$U:\R^2 \to \baR$ is a translation-invariant, smoothly approximable
standard potential. We choose $K$, $\psi$, $\bU$ and $u$ according 
to Definition \ref{defGapprox}. W.l.o.g.  we may assume $0 \in K$, 
$\bU = U$ and $u = 0$ on $\Ko$. We then let  $\ep >0$ so small that 
\begin{equation} \label{dec}
\ex  := \leb(\Kep \weg \KU) + \int_{\Ko^c} \tiu(x) dx  < \frac{1}{z\xi}.
\end{equation}
In addition to the function $\fc$ and the constants $\cK$ and $\cf$ 
introduced in Section \ref{constantshc} we also define
\begin{equation} \label{bpsi}
\cu  \, :=  \, \int_{\Ko^c} \tiu (x) |x|^2 dx  \quad \text{ and } \quad  
\cpsi \, := \, \|\psi\|\, \vee \int dx  \, \psi(x)(|x|^2 \vee 1).
\end{equation}
These constants are finite by our assumptions. Finally, we fix a Gibbs 
measure $\mu \in \G_{\XX}(U,z)$, a cylinder event $D \in \F_{\X,\La_{n'-1}}$ 
where $n'\in \N$, a translation distance  $\tau  \in [0,1/2]$, 
the translation direction $e_1$ and a real $\de > 0$. 

%=================================================================

\subsection{Decomposition of $\mu$ and the bond process}

For $n \in \N$ and $X \in \X$ we consider the bond set 
\[
\en(X)  \, := \,  E_{\Lan}(X) \, = \, 
\{ x_1x_2 \in E(X):  x_1 x_2 \cap \Lan \ne \emptyset\}.
\]
On $(\E_{\en(X)},\B_{\en(X)})$ we introduce the Bernoulli measure 
$\pin(.|X)$ with bond probabilities 
\[
(\tiu(b))_{b \in \en(X)} \quad \text{ where } \quad  \tiu(b) \, := \,  1-e^{-u(b)}, 
\]
using the shorthand notation $u(x_1x_2) := u(x_1-x_2)$ for $x_1,x_2 \in \R^2$. 
We  note that   $0 \le \tiu(b) < 1$ for all $b \in \en(X)$ as 
$0 \le u < \infty$. As remarked earlier $\pin(.|X)$ can be extended to 
a probability measure on $(\E, \F_{\E})$. For all $D \in \F_{\E}$ $~\pin(D|.)$ 
is $\F_{\X}$-measurable, so $\pin$ is a probability kernel from $(\X,\F_\X)$ 
to $(\E,\F_\E)$. 
\begin{lem} \label{umord}
Let $n \in \N$. We have  
\[   
\mu \otimes \nu_{\Lan}(\Gnb) = 1 \; \text{ and } \; \mu(\Gnb) = 1 
\; \text{ for } \; 
\Gnb  :=  \{X \in \XX: \sum_{b \in \en(X)} \tiu(b)   < \infty \}.
\]
\end{lem} 
For $X \in \Gnb$ the Borel-Cantelli lemma implies that every bond set  
is finite $\pin(.|X)$-a.s., so 
\[
\sideset{}{'}\sum_{B \subset \en(X)} \pin(\{B\}|X) \, = \, 1, 
\]
where the summation symbol $\sum'$  indicates that the sum
extends over finite subsets only. We have  
\[
\pin(\{B\}|X) \, 
= \, \prod_{b \in B} \tiu(b)  \prod_{b \in \en(X) \weg B} (1 - \tiu(b)) \,
= \, e^{-H^{{u}}_{\Lan}(X)}  \prod_{b \in B} (e^{u(b)}-1),  
\]
so for every $X \in \Gnb$ the Hamiltonian $H^u_{\Lan}(X)$ is finite, 
and thus the decomposition of the potential gives a corresponding
decomposition of the Hamiltonian 
\[
H^U_{\Lan}(X) \, = \,  H^{\bar{U}}_{\Lan}(X) - 
H^{{u}}_{\Lan}(X).
\]
Using \eqref{gibbsaequ} we conclude that for every 
$\F_\X \otimes \F_{\E}$-measurable function $f \ge 0$  
\begin{equation}  \label{desint} 
\begin{split} 
\int d\mu \otimes \pin \, f \quad  
&= \quad  \int \mu(d\bX) \frac 1  {Z_{\Lan}(\bar{X})} \int \nu_{\Lan}(dX|\bX) 
\sideset{}{'}\sum_{B \subset \en(X)} f(X,B) \\
&\hspace{ 3 cm} \times z^{\#X_{\Lan}} e^{-H^{\bar{U}}_{\Lan}(X)} 
 \prod_{b \in B} (e^{u(b)}-1). 
\end{split} \end{equation}
Here by  Lemma \ref{umord} on both sides we have  $X \in \Gnb$ with 
probability one, thus the equality follows from the above decomposition. 
If $f$ does not depend on $B$ at all, 
the integral on the left hand side of \eqref{desint} is just the $\mu$-expectation 
of $f$, as $\pin$  is a probability kernel, and from the right hand side
we learn that the perturbation $u$ of the smooth potential $\bU$ 
can be encoded in a bond process $B$ such that the perturbation affects
only those pairs of particles with $x_1x_2 \in B$.\\   

On $(\E_{\en(X)},\B_{\en(X)})$ we denote  the counting measure concentrated 
on finite bond sets by $\pinh(.|X)$. Again $\pinh$ can be considered as a 
probability kernel from $(\X,\F_\X)$ to $(\E,\F_\E)$. 
For all $\F_{\E}$-measurable functions $f \ge 0$ we have 
\[
\int \pinh(dB|X) f(B) \, = \, \sideset{}{'}\sum_{B \subset \en(X)} f(B).
\] 

%==========================================================================

\subsection{Generalised translation} \label{tgt} 

First of all, we need to augment each bond set $B$ by additional bonds 
between all particles that are close to each other. That is, for $n > n'$, 
$X \in \X$ and $B \subset \en(X)$  we introduce the $\Kep$-enlargement of $B$ by  
\[
B_+ := B \cup \{ x_1x_2 \in \en(X): x_1 - x_2 \in \Kep \}.
\]
We then consider the range of the $B_+$-cluster of $\La \in \Bo^2_b$
\begin{displaymath} 
\rnxb(\La) \; = \; \sup\{ |x'| : x' \in C_{X,B_+}(\La) \}. 
\end{displaymath}
\begin{lem} \label{lereich}
We have $ \quad \sup \limits_{n > n'} \int \mu \otimes \pin (dX,dB) 
   \, \rnxb(\La_{n'}) \,  < \, \infty$.
\end{lem}
By the Chebyshev inequality we therefore can choose an 
integer $R > n'$, such that for every $n >n'$ the event 
\[
\Gna \, := \, \{(X,B) \in \X \times \E: \rnxb(\La_{n'}) <  R, 
B \subset \en(X) \text{ finite} \} \, 
\in \, \F_{\X} \otimes \F_{\E}
\]
\[
\text{ has probability } \quad  
\mu \otimes \pin (\Gna) \, \ge \, 1- \de/2.
\]
For $n > R$ we define the functions
$q$, $Q$, $r$ and $\tn$ exactly as in Section \ref{tgthc}. 
For $X \in \X$, $B \in \en(X)$ and $x \in X$ we define 
$\anxb(x)$ to be a point of $C_{X,B_+}(y)$ such that 
\[
|\anxb(x)| \; \ge \; |x|, \quad 
\tn(|\anxb(x)|) \; 
= \; \min \{ \tn(|x'|) : x' \in C_{X,B_+}(y) \}
\]
and $\anxb(x)$ is a measurable function of $x$, $X$ and $B$.

%=========================================================

\subsection{Good configurations}

In order to deal with the hard core and the perturbation encoded in the bond 
process, we will introduce a transformation  
\[
\Tn: \X \times \E \to \X \times \E
\]
which is required to have the following properties: 
\begin{enumerate}
\item [(1)] Whenever $B$ is a set of bonds between particles in $X$, 
the transformed configuration  $(\tX,\tB) = \Tn(X,B)$  is constructed by 
translating every particle $x \in X$ by a certain distance  $\tnxb(x)$ 
in direction  $\einh$, and by translating bonds along with the corresponding 
particles. 
\item [(2)] Particles in the inner region $\La_{n'-1}$ are translated by $\tau \einh$,
and particles in the outer region  ${\Lan}^c$ are not translated at all.
\item [(3)] Particles connected by a bond in  $B$ are translated the same 
distance.  
\item [(4)] $\Tn$ is bijective, and the density of the transformed process with 
respect to the untransformed process under the measure $\nu \otimes \pinh$ 
can be calculated explicitly. 
\item [(5)] We have suitable estimates on this density and on 
$H^{\bU}_{\Lan}(\tX)-H^{\bU}_{\Lan}(X)$. For the last assumption we need particles 
within hard core distance to remain within hard core distance and particles at 
larger distance to remain at larger distance. 
\end{enumerate} 
Property (2) implies that the translation of the chosen cylinder event $D$ 
is the same as the transformation of $D$ by $\Tn$. Properties (3)-(5) are
chosen with a view to the right hand side of \eqref{desint}: If $\Tn$ has
these properties then the density of the transformed process with 
respect to the untransformed process under the measure $\mu \otimes \pin$
can be estimated. We will content ourselves 
with a transformation satisfying the above properties only for  
$(X,B)$ from a set of good configurations 
\begin{equation} \label{good}
 \Gn  \, :=  \, \big\{ (X,B) \in \Gna: 
  \sum \limits_{i=1}^5 \Sigma_i(n,X,B) <   1/2  \big\} \, 
\in \, \F_{\X} \otimes \F_{\E}.
\end{equation} 
The functions $\Sigma_i(n,X,B)$ will be defined whenever we want
good configurations to have a certain property. In Lemma
\ref{lebadsmall} we then will prove that the set of good configurations 
$\Gn$ has probability close to $1$ when $n$ is big enough.  
Up to that point we consider a fixed $n \ge R+1$.

%===================================================================

\subsection{Modifying the generalised translation}

The construction of the deformed translation $\Tn$ will go along 
the same lines as the corresponding construction in section \ref{secmodgen}. 
However, here we also have to consider bonds between particles, 
and by property (3) from the last section we know that we have
to translate not just particles, but whole  $B$-clusters.\\ 
For a rigorous recursive definition of $\Tn(X,B)$ we first consider the case
that $B$ is a finite subset of $\en(X)$. Let $\tnxb^0 := \tau_n(|.|)$,
$\cnxb^0$ the $B$-cluster of the outer region $\Lan^c$, $m=m(X,B)$ the number 
of different $B$-clusters of $X \weg \cnxb^0$ and $\taunxb^0 := 0$.
In the k-th construction step $(1 \le k \le m)$ let 
\[
\tnxb^{k} \, :=  \, \tnxb^{k-1} \wedge 
   \bigwedge_{x \in \cnxb^{k-1}}m_{x,\taunxb^{k-1}} \,
= \quad  \tnxb^{0} \wedge \bigwedge_{0 \le i < k}
   \bigwedge_{x \in \cnxb^{i}} m_{x,\taunxb^{i}},
\]
where the auxiliary function  $\mxt: \R^2 \to \bar{\R}$ is defined as
in Section \ref{secmodgen}.  
Let the pivotal point $\pnxb^k$ be the point of 
$X \weg (\cnxb^0 \cup \ldots \cup \cnxb^{k-1})$ 
at which the minimum of $\tnxb^{k}$ is attained. 
If there is more than one such point then take the smallest point 
with respect to the lexicographic order for the sake of definiteness. 
Let  $\taunxb^k := \tnxb^{k}(\pnxb^k)$  be the corresponding minimal value of $\tnxb^k$,
$\cnxb^k$ the $B$-cluster of the point $\pnxb^k$ and  $\Tnxb^k := id + \tnxb^k \einh$.
For $k = m+1$ we can still define  $\tnxb^{m+1}$, but then the recursions stops as 
$X \weg (\cnxb^0 \cup \ldots  \cup \cnxb^{m}) = \emptyset$.
In the above notations we will omit dependence on $X$ and $B$ if it is 
clear which configuration is considered.
Now  for $x \in \cnxb^k$ let  $\tnxb(x) := \taunxb^{k}$ be 
the deformed translation distance function and let 
\[ 
\begin{split}
\Tnb(X) \, &:= \, \bigcup_{k=0}^{m(X,B)} (\cnxb^k + \taunxb^k \einh) =
\{x + \tnxb(x) \einh : x \in X \}  \quad \text{ and }\\
\Tnx(B) \, &:= \, \{ (x + \tnxb(x)\einh)(x' + \tnxb(x')\einh): xx' \in B \}. 
\end{split}
\]
If $B$ is not a finite subset of $\en(X)$ we define $\Tnb = id$ and $\Tnx = id$.
The deformed transformation can now be defined to be 
\[
\Tn: \X \times \E \to \X \times \E, \quad \Tn(X,B) \, := \, (\Tnb(X),\Tnx(B)).  
\]
Using Lemma \ref{lemeas} one can show that all above objects 
are measurable with respect to the considered $\si$-algebras. 
In the rest of this section we will convince ourselves that the above 
construction has indeed the required properties. 
\begin{lem} \label{leinnenaussen}
For good configurations  $(X,B) \in \Gn$ we have    
\begin{equation} \label{innenaussen}
(\Tnb X - \tau \einh)_{\La_{n'-1}} \, = \, X_{\La_{n'-1}}\quad \text{ and }
\quad (\Tnb X)_{{\Lan}^c} \, = \, X_{{\Lan}^c}.
\end{equation}
\end{lem}
\begin{lem} \label{lebij}
The transformation $\Tn: \X \times \E \to \X \times \E$ is bijective.
\end{lem}
In the proof of Lemma \ref{lebij} we again construct the inverse 
of $\Tn$, which is needed in the proof of the following lemma.
There we will also show that
\begin{equation} \label{densdef}
\ph (X,B) \, := \, 
\prod_{k = 1}^{m(X,B)} \big|1 +  \partial_{1} \tnxb^{k} (\pnxb^k)\big| 
\end{equation}
is well defined $\nu_{\Lan} \otimes \pinh (\,. \, |\bar{X})$-a.s., in that 
the considered derivatives exist. 
\begin{lem} \label{ledensp}
For every  $\bar{X} \in \X$ and every $\F_{\X}  \otimes \F_{\E}$-measurable 
function  $f \ge 0$ 
\begin{equation} \label{dens}
\int d\nu_{\Lan}\otimes \pinh(.|\bar{X}) \, (f \circ \Tn \cdot \ph) \,
= \, \int d\nu_{\Lan}\otimes \pinh(.|\bar{X}) \, f.
\end{equation}
\end{lem}
We also need  the backwards translation. So let
$\iTn$, $\iTnB$, $\iTnx$ and $\bph$ be defined analogously to the above objects, 
where now $\einh$ is replaced by $-\einh$. The previous lemmas apply analogously 
to this deformed backwards translation. We note that $\iTn$ is not the inverse 
of $\Tn$.

%=========================================================================

\subsection{Final steps of the proof}  \label{Pott}

From \eqref{desint} and Lemma \ref{ledensp} we deduce  
\[ \begin{split}
\mu &\otimes \pin (\Tn (D \cap \Gn))\\ 
&= \,  \int \mu(d\bX) \frac 1 {Z_{\Lan}(\bar{X})} 
  \int \nu_{\Lan} \otimes \pinh(dX,dB|\bX)\\ 
&\hspace{0.3 cm} 1_{\Tn(D \cap \Gn)} \circ \Tn(X,B) \, z^{ \# (\Tnb X)_{\Lan}} 
  \ph(X,B) e^{-H^{\bar{U}}_{\Lan}(\Tnb X)} \prod_{b \in \Tnx B} (e^{u(b)}-1).
\end{split} \]
Here we have identified $D$ and $D \times \E$. By the bijectivity of $\Tn$ 
from Lemma \ref{lebij},  by \eqref{innenaussen} and by  construction of $\Tnx$ 
the above integrand simplifies to 
\[
 1_{D \cap \Gn} (X,B)  \, z^{\# X_{\Lan}} \, 
  e^{\log \ph (X,B) - H^{\bU}_{\Lan}(\Tnb X)} \prod_{b \in B} (e^{u(b)}-1).
\]
The backwards transformation  $\iTn$ can be treated analogously, 
hence
\[ \begin{split}
&\mu \otimes \pin (\iTn (D \cap \Gn))\, 
 + \,  \mu \otimes \pin (\Tn (D \cap \Gn)) \, 
 - \, \mu \otimes \pin(D \cap \Gn) \\ 
&= \,  \int \mu(d\bX) \frac{1}{Z_{\La_{n}}(\bar{X})} 
  \int \nu_{\Lan} \otimes \pinh(dX,dB|\bX) \, 1_{D \cap \Gn} (X,B)  
 \, z^{\# X_{\Lan}}  \prod_{b \in B} (e^{u(b)}-1) \\  
& \hspace{2 cm} \times  \Big[
  e^{\log \bph(X,B) -H^{\bU}_{\La_{n}}(\bar{T}_{n,B} X)} 
 + e^{\log \ph(X,B) -H^{\bU}_{\Lan}({T}_{n,B} X)}  - e^{ -H^{\bU}_{\Lan}(X)} \Big].
\end{split} \]
We note that  for  $(X,B) \in \Gn$ we have 
\[
\begin{split}
&e^{\log \bph(X,B) -H^{\bU}_{\La_{n}}(\bar{T}_{n,B} X)} 
 + e^{\log \ph(X,B) -H^{\bU}_{\Lan}({T}_{n,B} X)} \\ 
& \qquad \ge \, 2 \, e^{\frac{1}{2}( \log \bph(X,B) + \log \ph(X,B)  
 -  H^{\bU}_{\La_{n}}(\bar{T}_{n,B} X) -  H^{\bU}_{\Lan}({T}_{n,B} X))}\\
& \qquad \ge \, 2 \, e^{- \frac 1 2 -  H^{\bU}_{\Lan}(X)} \quad 
\ge \quad   e^{-  H^{\bU}_{\Lan}(X)},  
\end{split}
\]
where we have used the convexity of the exponential function in the first step
and the following estimates in the second step:  
\begin{lem} \label{ngross}
For $(X,B) \in \Gn$ we have     
\begin{equation} \label{taylor}
H^{\bU}_{\La_{n}}(\iTnB X) + H^{\bU}_{\La_{n}}(\Tnb X) \,
\le \, 2 H^{\bU}_{\Lan}(X) + 1/2 \quad \text{ and }
 \end{equation}
\begin{equation} \label{dichtenklein}
\log \bph(X,B)  + \log \ph(X,B) \, \ge \, - 1/2. 
\end{equation}
\end{lem}
Hence we have shown that 
\begin{equation} \label{hingeorgii}
\mu \otimes \pin (\iTn (D \cap \Gn)) \, + \, 
\mu \otimes \pin (\Tn (D \cap \Gn)) \, 
\ge \, \mu \otimes \pin (D \cap \Gn) .
\end{equation}
In \eqref{hingeorgii} we would like to replace $D \cap \Gn$ by $D$, 
and for this we need $\Gn$ to have high probability:
\begin{lem} \label{lebadsmall}
If $n \ge R + 1$ is chosen big enough, then $\mu \otimes \pin (\Gn^c) \le \de$.
\end{lem}
For the proof of Theorem \ref{sym} we choose such an $n \ge R+1$. The 
rest of the argument is then the same as that in Section \ref{Potthc}.

% %=======================================================================
% %
% %======================================================================

\section{Proof of the lemmas from Section \ref{secproofsym}} 
\label{secleproofsym}

\subsection{Convergence of energy sums: Lemma \ref{umord}}

Let $n \in \N$. For every $X \in \X$ we have 
\begin{displaymath}
H^{\tiu}_{\Lan}(X) \, = \, \sum_{b \in \en(X)} \tiu(b) \, 
\le  \, \sideset{}{^{\neq}}\sum_{x_1, x_2 \in X} 
1_{\{x_1 \in \Lan\}} \, \tiu(x_1-x_2), \quad \text{ and so }  
\end{displaymath}
\begin{displaymath}
\int \nu_{\Lan}(dX|\bX) H^{\tiu}_{\Lan}(X) \, 
\le  \, \int_{\Lan} dx_1 \Big( \int_{\Lan} dx_2  \tiu(x_1-x_2)
+  \sum_{x_2 \in \bX_{\Lan^{\;c}}} \tiu(x_1-x_2) \Big)
\end{displaymath}
for all  $\bX \in \X$. By Lemma  \ref{lekorab} we get
\begin{displaymath}
\begin{split}
\int \mu \otimes \nu_{\Lan}(dX) &H^{\tiu}_{\Lan}(X) 
\le  \int_{\Lan} dx_1 \Big( \int_{\Lan} dx_2  \tiu(x_1-x_2)
+  z \xi \int_{\Lan^{\;c}} dx_2  \tiu(x_1-x_2) \Big)\\
&\le \, \int_{\Lan} dx_1 (1 + z\xi) \ex  \quad 
\le \quad 4n^2  (1 + z\xi) \ex \quad < \quad \infty, 
\end{split}
\end{displaymath}
where we have estimated the integrals over  $x_2$ by $\ex$ using \eqref{dec}. 
Thus we have proved the first assertion. However, $\mu$ is absolutely 
continuous with respect to $\mu \otimes \nu_{\Lan}$, which follows from  
\eqref{gibbsaequ} and the definition of the conditional Gibbs distribution. 
So the first assertion implies the second one.

%================================================================

\subsection{Cluster bounds: Lemma \ref{lereich}}

Let us refine the argument of Section \ref{seclereichhc} as follows. 
For $n>n'$, $X \in \X$ and $B \subset \en(X)$ we consider a path 
$x_0,...,x_m$ in the graph $(X,B_+)$ such that $x_0 \in \La_{n'}$, 
and we consider an integer $k$ such that $1 \le k \le m$ and the bond 
$x_{k-1}x_k$ has maximal $|.|$-length among all bonds on the path. We have  
\[ 
|x_m| \, \le \,  |x_0| + \sum_{i=1}^m |x_i-x_{i-1}| \,  
\le \, n' + m   |x_k-x_{k-1}|. 
\] 
By considering all paths and bonds of maximal length we obtain
\[
\rnxb (\La_{n'}) \, \le \, n' +  \sum_{m \ge 1} \sum_{k = 1}^m 
 \quad \sideset{}{^{\neq}} \sum_{x_0,\ldots ,x_m \in X}   1_{\{x_0 \in \La_{n'} \}}  
 m |x_k-x_{k-1}| \prod_{i=1}^m   1_{ \{x_i x_{i-1} \in B_+\} }.
\]
Under the Bernoulli measure  $\pin(dB|X)$, the events 
$\{x_i x_{i-1} \in B_+ \}$ are independent, and for 
$g := 1_{\Kep \weg \KU} +  \tiu$ we have 
\begin{equation} \label{pig}
\int \pin(dB|X)  1_{ \{x_i x_{i-1} \in B_+\} }  \, 
\le \,1_{\KU}(x_i-x_{i-1}) +g(x_i-x_{i-1}). 
\end{equation} 
Using the hard core property  \eqref{hardcore}
and Lemma \ref{lekorab} we thus find
\[ \begin{split}
R_n \, &:= \, \int \mu(dX) \int \pin(dB|X) \rnxb (\La_{n'}) \, - \, n'\\
&\le \, \sum_{m \ge 1} \sum_{k = 1}^m \int \mu(dX) \, \sideset{}{^{\neq}}
 \sum_{x_0,\ldots ,x_m \in X} 1_{\{x_0 \in \La_{n'} \}}
  m |x_k-x_{k-1}| \prod_{i=1}^m g(x_i-x_{i-1})\\
&\le \, \sum_{m \ge 1} \sum_{k = 1}^m \, (z\xi)^{m+1} \, 
   \int dx_0 \ldots  dx_{m}   1_{\{x_0 \in \La_{n'} \}}
  m |x_k-x_{k-1}| \prod_{i=1}^m g(x_i-x_{i-1}).
\end{split} \] 
Setting $c_g := (1 + \cK^2)\ex  + \cu$ we conclude from \eqref{dec} and \eqref{bpsi} 
that
\begin{equation} \label{intg}
\int g(x)|x|\, dx \, \le \,  \int g(x)(1 + |x|^2)\, dx \, \le \, c_g
 \quad \text{ and } \quad 
\int g(x) \, dx \, \le \, \ex, 
\end{equation}
hence we can estimate the integrals over $dx_i$ in the above expression 
beginning with  $i =m$. This gives $m-1$ times a factor 
$\ex$ and once a factor $c_g$. Finally the integration over  $dx_0$ 
gives an additional factor  $\la^2(\La_{n'}) = (2n')^2$. Thus   
\[
  R_n \, \le  \, (2n'z\xi)^2 c_g \sum_{m \ge 1}  
  m^2 (\ex z\xi )^{m-1}.
\]
The last sum is finite  because $\ex z\xi < 1$.

%=====================================================================

\subsection{Properties of the deformed translation: Lemma~\ref{leinnenaussen}}

We will show properties of the construction which are analogous
to properties of the corresponding objects from the proof of the special case
in Sections \ref{secleconshc} and \ref{secpdt}. Additionally
we need  a way to calculate the translation distance of an 
arbitrary particle $x \in \ecn^k$ without knowing $\epn^k$. This can be done
using the first relation of the following lemma. 
\begin{lem} \label{letrd}
For $X \in \X$, finite $B \subset E(X)$, $k \ge 0$, $x,x' \in X$ and 
$s \in [-1,1]$
\begin{align}
&\etaun^k \, = \, \etn^{k+1}(x) \quad \text{ if } x \in \ecn^k, \label{ktk}\\ 
&\etaun^k \le \etaun^{k+1}, \label{mono}\\ 
&\etn^k \text{ is 1/2-$\einh$-Lipschitz-continuous and piecewise  
cont. $\einh$-differentiable}, \label{lip} \\  
&\eTn^k \text{ is $\lear$-increasing and bijective}, \label{wabi}\\ 
&x-x' \in \Ko  \quad 
\Rightarrow \quad \tnxb(x) \,= \,\tnxb(x'), \label{vercon}\\
&x-x' \notin \Ko \quad \Rightarrow \quad 
x-x'+ s(\tnxb(x)-\tnxb(x'))\einh) \notin \Ko, \label{greps}\\
&\etaun^k \, \le \, \etn^0(x) \quad 
\text{ for all } x \in \ecn^{k'} \text{ such that } k' \ge k 
, \label{tauoben}\\
&\etaun^k \, \ge \, \etn^0(\anxb(\epn^k))    \quad \text{ if } 
  (X,B) \in \Gn,  \label{tauunten}
\end{align}
\end{lem}
\Bew
For \eqref{ktk} let $x \in \ecn^k$. By definition of  $\epn^k$ we have 
$\etn^{k}(x) \ge \etaun^k$, so 
\[ 
\etn^{k+1}(x) \, =\, \etn^{k}(x) \wedge 
 \bigwedge_{ x' \in \ecn^k} m_{x',\etaun^k}(x) \, = 
\, \etaun^k,
\]
where we have also used $m_{x',\etaun^k}(x) \ge \etaun^k$ 
and $m_{x,\etaun^k}(x) = \etaun^k$. The other assertions can be shown 
as in Sections \ref{secleconshc} and \ref{secpdt}. Here for the proof of 
\eqref{vercon} and \eqref{greps} we have to use  \eqref{ktk}, and  the key 
observations are the following: For $i \le j$, $x_i \in \ecn^i$ and  
$T^{j+1}_{n,s} := id + s \cdot \etn^{j+1} \einh$ we have  
\[ \begin{split} 
&\alle x \in \Ko(x_i) \cap \ela^i: \quad  
\etn^{j+1}(x) \, =\, \etn^{i}(x) \wedge \bigwedge_{i \le k \le  j} \, 
 \bigwedge_{ x' \in \ecn^k} m_{x',\etaun^k}(x) \, =  \, \etaun^i,\\
&T^{j+1}_{n,s}(\ela^i) \, = \, \ela^i + s \etaun^i \einh \quad \text{ and } 
\quad T^{j+1}_{n,s}(\ela^i \weg \Ko(x_i)) \, = \, \ela^i  \weg \Ko(x_i)  
+ s \etaun^i \einh.  
\end{split}
\]  
To obtain  \eqref{tauunten} here we specify the function  
\begin{equation} \label{sig1} 
\Sigma_1(n,X,B)\, :=  \sum_{x,x' \in X} 1_{\{|x| \le |x'|\}} 
 1_{\{ x \stackrel{X,B_+}{\longleftrightarrow} x'\}} 4
 \big(\tn(|x|-\cK) - \tn(|x'|)\big)^2 \cf^2
\end{equation}
used in the definition of $\Gn$. \qed\\

\noindent 
Lemma  \ref{leinnenaussen} follows from  \eqref{tauoben} and \eqref{tauunten}, 
 just as in the proof of Lemma~\ref{leinnenaussenhc}.

%======================================================================

\subsection{Bijectivity of  the transformation: Lemma \ref{lebij}}
\label{secbij}

The construction of the inverse transformation is analogous to the one in 
Section~\ref{secbijhc}. Let $\tX \in \X$ and $\tB \subset \en(\tX)$ be finite.
Let $\tcnxb^0$ be the $\tB$-cluster of $\Lan^c$, $\tilde{m}=\tilde{m}(\tX,\tB)$ 
the number of different $\tB$-clusters of $\tX \weg \tcnxb^0$,  
$\ttnxb^0 = \tau_n(|.|)$  and $\ttaunxb^0 := 0$. 
In the k-th construction step $(k \ge 1)$ let
\[
\ttnxb^{k} \, :=  \, \ttnxb^{k-1} \wedge 
   \bigwedge_{x \in \tcnxb^{k-1}-\ttaunxb^{k-1}} m_{x,\ttaunxb^{k-1}}.
\]
Let $\tTnxb^k = id + \ttnxb^k \einh$ and $\tpnxb^k$ be the point of 
$\tX \weg (\tcnxb^0 \cup \ldots  \cup \tcnxb^{k-1})$ at which the minimum
of $\ttnxb^k \circ (\tTnxb^k)^{-1}$ is attained. If there is more than one 
such point then take the point $x$ such that $(\tTnxb^k)^{-1}(x)$ is minimal 
with respect to the lexicographic order $\le$. Let  
$\ttaunxb^k := \ttnxb^k \circ (\tTnxb^k)^{-1}(\tpnxb^k)$ be  the corresponding 
minimal value and $\tcnxb^k$ be the $\tB$-cluster of the pivotal point $\tpnxb^k$. 
The recursion stops for  $k = \tilde{m} + 1$.
In the above notations we will omit dependence on $\tX$ and $\tB$ 
if it is clear which configuration is considered. 
We need to show that the above construction is well defined, i.e. 
that   $\tTnxb^k$ is invertible in every step. Furthermore we 
need some more properties of the construction. All this is done in the
following lemma: 
\begin{lem} \label{leeiginv}
Let $\tX \in \X$, $\tB \subset \en(\tX)$ finite and $k \ge 0$. Then
\begin{align}
& \ettn^k \text{ is 1/2-$\einh$-Lipschitz-continuous}, \;
\etTn^k \text{ is bijective and $\lear$-increasing}, \label{lipinv}\\
& (\etTn^k)^{-1} + \ettn^k \circ (\etTn^k)^{-1} \einh  \, = \,id, 
 \label{TSinv}\\
& \alle c \in \R ,x \in \R^2: \, \ettn^k \circ (\etTn^k)^{-1}(x) 
 \,\ge\,c\, 
  \Leftrightarrow \, \ettn^k(x -  c\einh) \, \ge \, c, \label{tsinv}\\ 
& \ettn^{k} \, \le \, \ettn^{k-1} \quad \text{ and } \quad  
 \ettaun^{k-1} \, \le \, \ettaun^{k}, \label{monoinv}\\
& \alle x \in \etcn^{k}: \, \ettn^{k+1} \circ(\etTn^{k+1})^{-1}(x) \, 
= \, \ettaun^k.
\label{kek}
\end{align} 
\end{lem}
\Bew
Assertions \eqref{lipinv} - \eqref{monoinv} can be shown exactly as 
the corresponding assertions from Lemma \ref{leeiginvhc}. 
For \eqref{kek} let  $x \in \etcn^k$. We have
\[ \begin{split}
\ettn^{k} &\circ(\etTn^{k})^{-1}(x) \,\ge \,\ettaun^k 
\quad \Rightarrow \quad 
\ettn^{k}(x-\ettaun^k \einh ) \,\ge \,\ettaun^k \\
&\Rightarrow \quad \ettn^{k+1}(x - \ettaun^k \einh) \,= \, \ettaun^k 
\quad \Rightarrow \quad 
\ettn^{k+1} \circ(\etTn^{k+1})^{-1}(x) \,= \, \ettaun^k,
\end{split} \]
where the first statement holds by definition, and the implications 
follow from \eqref{tsinv}, 
$x - \ettaun^k \einh \in \etcn^k - \ettaun^k \einh$ and 
\eqref{TSinv} respectively. \qed\\ 

\begin{figure}[!htb] 
\begin{center}
\psfrag{1}{$\epn^0$}
\psfrag{2}{$\etpn^0$}
\psfrag{5}{$\epn^2$}
\psfrag{6}{$\etpn^2$}
\psfrag{7}{$\epn^1$}
\psfrag{8}{$\etpn^1$}
\psfrag{a}{$\ecn^0$}
\psfrag{b}{$\etcn^0$}
\psfrag{e}{$\ecn^2$}
\psfrag{f}{$\etcn^2$}
\psfrag{g}{$\ecn^1$}
\psfrag{h}{$\etcn^1$}
\psfrag{A}{$\etaun^0$}
\psfrag{B}{$\ettaun^0$}
\psfrag{E}{$\etaun^2$}
\psfrag{F}{$\ettaun^2$}
\psfrag{G}{$\etaun^1$}
\psfrag{H}{$\ettaun^1$}
\psfrag{x}{$\Tn: (X,B) \mapsto (\tX,\tB)$}
\psfrag{y}{$\tTn: (\tX,\tB) \mapsto (X,B)$}
\includegraphics[scale=0.5]{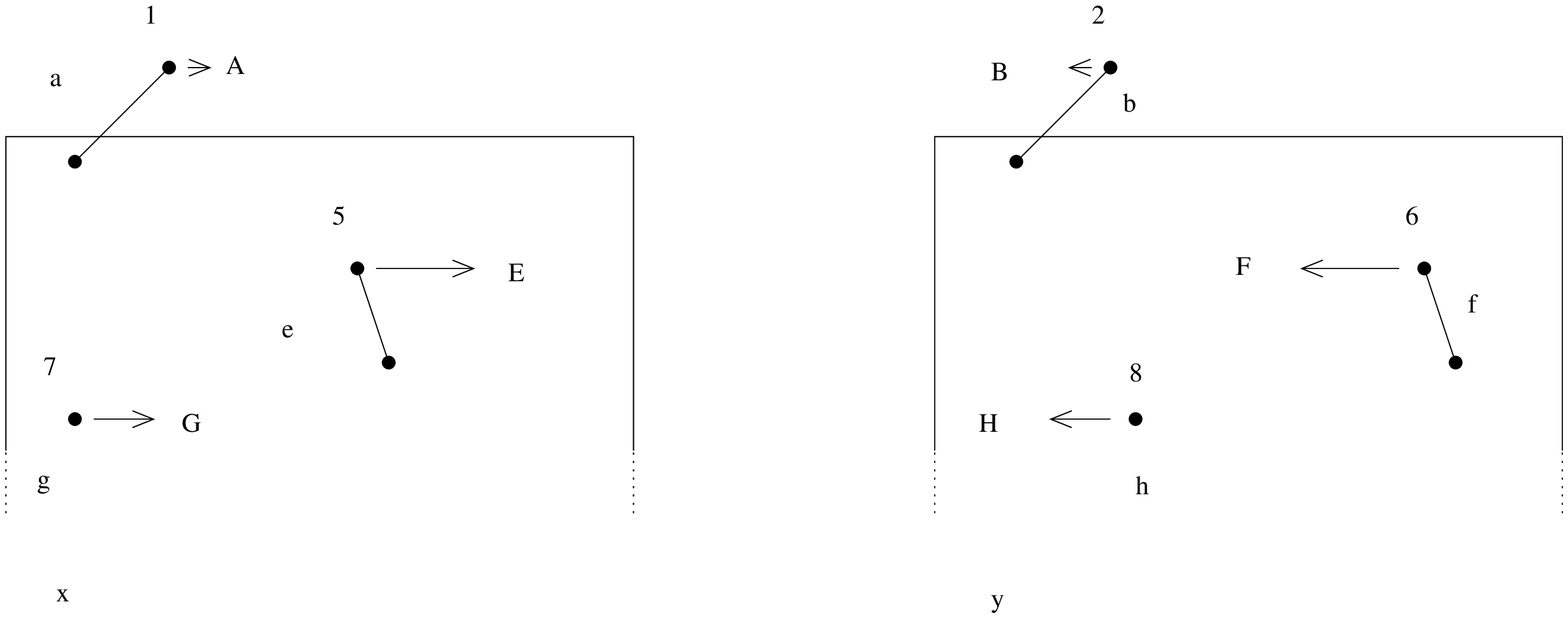}
\end{center}
\caption{Construction of the inverse  $\tTn$ of $\Tn$.} \label{figinv}
\end{figure}

\noindent
For $x \in \tcnxb^k$ let $\ttnxb(x) := \ttaunxb^{k}$ be the distance
the particle $x$ is translated. We define 
\[ 
\begin{split}
\tTnb(\tX) \, &:= \, \bigcup_{k=0}^m (\tcnxb^k - \ttaunxb^k \einh) =
\{ x - \ttnxb(x) \einh : x \in \tX \}  \quad \text{ and }\\ 
\tTnx(\tB) \, &:= \, \{ (x - \ttnxb(x)\einh)(x' - \ttnxb(x')\einh): xx' \in \tB \}. 
\end{split}
\]
Now if $\tB$ is a not a finite subset of  $\en(\tX)$ we define   
$\tTnb = id$ and $\tTnx = id$. $\Tn$ is then defined by 
\[
\tTn: \X \times \E \to \X \times \E, \quad \tTn(X,B) \, := \, (\tTnb(X),\tTnx(B)).  
\]
By Lemma \ref{lemeas} we see again that all above objects are 
measurable with respect to the considered $\si$-algebras. The following 
two lemmas are the key to show that $\tTn$ is indeed the inverse of $\Tn$. 
The proofs differ from the proofs of Lemmas \ref{leminimhc} and 
\ref{leminiminvhc} only, in that we have to use \eqref{ktk} and \eqref{kek},
whenever we want to calculate the translation distance of a point 
explicitly. 
In the first lemma we consider $X \in \X$, finite $B \subset \en(X)$, $\etn^k$, 
$\eTn^k$, $\ecn^k$, $\epn^k$ and  $\etaun^k$ $~(0 \le k \le m)$ as in the 
construction of $\Tn(X,B)$, and we define $(\tX,\tB) := \Tn(X,B)$,  
$\etpn^k := \epn^k + \etaun^k \einh$ and 
$\etcn^k := \ecn^k + \etaun^k \einh$, see Figure \ref{figinv}.

\begin{lem} \label{leminim}
Let $1 \le k \le m$. For every   
$\tx \in \tX \weg (\etcn^0 \cup \ldots  \cup \etcn^{k-1})$ we have  
\[
\etn^k \circ (\eTn^k)^{-1}(\etpn^k) \,\le \,\etn^k \circ (\eTn^k)^{-1}(\tx).
\] 
For all $\tx$ for which equality occurs we have  
$(\eTn^k)^{-1} (\etpn^k) \le  (\eTn^k)^{-1}(\tx)$. 
\end{lem} 
For the second lemma we consider   $\tX \in \X$, finite $\tB \subset \en(\tX)$, 
$\ettn^k$, $\etTn^k$, $\etcn^k$, $\etpn^k$ and  $\ettaun^k$ 
$~(0 \le k \le \tilde{m})$ as in the construction of $\tTn(\tX,\tB)$, 
and we define $(X,B) := \tTn(\tX,\tB)$, $\epn^k := \etpn^k - \ettaun^k \einh$ and 
$\ecn^k := \etcn^k - \ettaun^k \einh$, see Figure~\ref{figinv}.
\begin{lem} \label{leminiminv}
Let $1 \le k \le m$. For every 
$x \in X \weg (\ecn^0 \cup \ldots  \cup \ecn^{k-1})$ we have 
\[
\ettn^k (\epn^k) \le \ettn^k (x).
\] 
For all $x$ for which equality occurs we have $\epn^k \le x$.  
\end{lem} 
Now the following lemma follows exactly as in the proof of Lemma \ref{leinvhc}. 
\begin{lem}\label{leinv}
On $\X \times \E$ we have $ \quad  \tTn \circ \Tn \, = \, id \quad $ 
and  $ \quad  \Tn \circ \tTn \, = \, id$.
\end{lem}

%==================================================================

\subsection{Density of the transformed process: Lemma \ref{ledensp}}

By definition the left hand side of \eqref{dens} equals 
\[ 
e^{-4 n^2} \sum_{k \ge 0} \frac{1}{k!} I(k), \quad
\text{ where }  \,
  I(k) \, = \, \int_{{\Lan}^{k}} dx
   \sideset{}{'}\sum_{B \subset \en(\bX_x)} (f \circ \Tn \cdot \ph) (\bX_x,B),
\]
using the shorthand notation $\bX_x = \{x_1,\ldots ,x_k\} \cup \bX_{\Lan^c}$.
We would like to fix the bond set  $B$ before we choose the positions 
$x_i$ of the particles. Thus we introduce bonds between indices of 
particles instead of bonds between particles. Let  $\Nk := \{1,\ldots ,k\}$,    
\[ 
\bX^k := \Nk \cup \bar{X}_{\Lan^c} \, 
\text{ and } \,  
E_n(\bX^k) := \{ x_1x_2 \in E(\bX^k): 
x_1x_2 \cap \Nk \neq \emptyset \}.
\]
For $B \subset E_n(\bX^k)$ and $x  \in \Lan^{\;I}$ 
$\,(I \subset \Nk)$ we define $B_x$ to be the bond set constructed from $B$ 
by replacing the point  $i \in I$ by $x_i$ in every bond of $B$ and by deleting 
every bond $B$ that contains a point $i \in \Nk \weg I$. Analogously let 
$\bX_x := \{ x_i: i \in I\} \cup \bX_{\Lan^c}$ be the configuration corresponding
to the sequence and let $(\bX,B)_x := (\bX_x,B_x)$. Using this notation we obtain
\[
I(k)  \,=  \sideset{}{'}\sum_{B \subset E_n(\bX^k)} I(k,B),   
\quad \text{ where } \,   I(k,B)  \,:= \,  \int_{{\Lan}^{k}} dx \,
   (f \circ \Tn \cdot \ph)(\bX,B)_x.   
\]
To compute $I(k,B)$ we need to calculate  $\Tn(\bX,B)_x$, and for this 
we must identify the points $P^i_{n,\bX_x,B_x}$ among the particles $x_j$. 
So let $m_B$ be the number of different $B$-clusters of 
$\bX^k \weg C_{\bX^k,B}(\Lan^c)$, 
$C_{\bX^k,B}(\per(0)) := C_{\bX^k,B}(\Lan^c)\cap\Nk$ 
and $\Pi(B)$ be the set of all mappings 
$\per:\{1,\ldots ,m_B\}\to (\bX^k \weg C_{\bX^k,B}(\Lan^c))$ such that 
every $\per(i)$ is in a different $B$-cluster. For $\per \in \Pi(B)$ 
let 
\[ \begin{split}
\akb \, &:= \, \big\{ x \in\Lan^{\;k}: \alle 1 \le j \le m_B: 
 x_{\per(j)} = P^j_{n,\bX_x,B_x}\big\}\quad \text{ and } \\
\takb \, &:= \, \big\{ x \in\Lan^{\;k}: \alle 1 \le j \le m_B:
 x_{\per(j)} = \tilde{P}^j_{n,\bX_x,B_x}\big\},
\end{split}\]
where  $\tilde{P}^j_{n,\bX_x,B_x}$ are the pivotal points from the construction 
of the inverse transformation in Subsection \ref{secbij}.  
Now we can write
\[
I(k,B) \, = \,\sum_{\per \in \Pi(B)} \int_{{\Lan}^{k}} dx \, 
        1_{\akb}(x) (f \circ \Tn \cdot \ph)(\bX,B)_x   
\]
and we denote the summands in the last term by  $I(k,B,\per)$. 
If  $x \in \akb$ we can derive a simple expression for $\Tn(\bX,B)_x$. 
For $x \in \Lan^{\;k}$ and $\per \in \Pi$ we define a formal transformation 
$T_{B,\per}(x)  :=  (T^i_{B,\per,x} (x_i))_{1 \le i \le k}$, 
where 
\[
t^{\per(j)}_{B,\per,x}  :=  t_{n,(\bX,B)_{x^{\per,j-1}}}^{j}, 
T^{\per(j)}_{B,\per,x}  :=  id + t^{\per(j)}_{B,\per,x}  \einh \;
\text{ and } \;   
T^i_{B,\per,x} :=  id + t^{\per(j)}_{B,\per,x}(x_{\per(j)})\einh
\]
for $0 \le j \le m_B$ and $i \in C_{\bX^k,B}(\per(j)), i \neq \per(j)$. Here 
$x^{\per,j}$ is defined to be the subsequence of $x$ corresponding to the index set 
$C_{\bX^k,B}^{\per,j} := \bigcup_{i \le j} C_{\bX^k,B}(\per(i))$.
Clearly, for $i \in C_{\bX^k,B}(\per(j))$,  $T^{i}_{B,\per,x}$ doesn't depend
on all components of $x$, but only on those $x_l$ such that 
$l \in C_{\bX^k,B}^{\per,j-1}$ and additionally on  $x_{\per(j)}$ if $i \neq \per(j)$. 
By definition we now have 
\begin{equation} \label{at}
x \in \akb \, \Rightarrow  \, \left\{ 
\begin{aligned} &\T_n(\bX,B)_x \, 
= \, (\bX, B)_{T_{B,\per}(x)} 
 \quad \text{ and }\\  
\, &T^j_{n,\bX_x,B_x} \, 
 = \, T^{\per(j)}_{B,\per,x}
\, \text{ for all } j \le m_B.
\end{aligned} \right. \end{equation}
Furthermore we observe that for all  $x \in (\R^2)^{k}$  we have 
\begin{equation} \label{trapo}  
x \in \akb \quad \Leftrightarrow \quad T_{B,\per}(x) \in \takb.
\end{equation}
This can be shown exactly as \eqref{trapohc} in Section \ref{secledensphc}.
Let $g: (\R^2)^{k} \to \R$, $g(x) :=  1_{\takb}(x) f(\bX_x,B_x)$.
Then \eqref{at} and \eqref{trapo} imply 
\[
I(k,B,\per) \, = \,  \Big[ \prod_{j=0}^{m_B} 
 \Big(   \prod_{i \in C_{\bX^k,B}(\per(j))} \int dx_i \Big)\,
 \big|1 + \partial_{1} t^{\per(j)}_{B,\per,x}(x_{\per(j)})\big| \Big] \, 
 g(T_{B,\per}(x)),  
\]
where we have also inserted the definition  of $\ph$ \eqref{densdef}.
Now we transform the integrals. For $j=m_B$ to $1$ and $i \in C_{\bX^k,B}(\per(j))$ 
we substitute $x_i' := T^i_{B,\per,x} x_i$. For $i \neq \per(j)$  $~T^i_{B,\per,x}$ 
is a translation by a constant vector, so $dx_i' = dx_i$. 
For $i = \per(j)$ the Lebesgue transformation theorem \eqref{lebtrsatz} gives 
\[
dx_{\per(j)}' \, = \, \big| 1 + \partial_{1} t^{\per(j)}_{B,\per,x}  
(x_{\per(j)})\big| dx_{\per(j)}
\] 
as in Section \ref{secledensphc}. Thus 
\[
I(k,B,\per) \, = \,   \Big(\prod_{j=0}^{m_B} 
   \prod_{i \in C_{\bX^k,B}(\per(j))} \int dx'_i 
  \Big) \, g(x') \, =  \, \int_{{\Lan }^k} dx \, 1_{\takb}(x) f(\bX_x,B_x)
\]
and we are done as the same arguments show that the right hand side of 
\eqref{dens} equals 
\[ 
 e^{-4 n^2} \sum_{k \ge 0} \frac{1}{k!}
\sideset{}{'}\sum_{B \subset E_n(\bX^k)} \sum_{\per \in \Pi(B)} 
 \int_{{\Lan }^k} dx \,   1_{\takb}(x) f(\bX_x,B_x).
\]
Combining the above ideas with the reasoning in Section \ref{secledensphc}
also shows  that the density function is well defined. 

%=======================================================================

\subsection{Key estimates: Lemma \ref{ngross}}

For all  $x \in \R^2$ and $\vartheta \in [-1,1]$ such that  
$x+ s\einh \notin \Ko$ for all $s \in [-\vartheta,\vartheta]$ we have
\[
\bU(x+\vartheta \einh) + \bU(x- \vartheta \einh) - 2 \bU(x) \, 
\le \,  \sup_{s \in [-\vartheta,\vartheta]} 
  \partial^2_{1} \bU(x+s\einh) \vartheta^2 \, \le\, \psi(x) \vartheta^2
\]
by Taylor expansion of  $\bU$ at $x$ using the $\einh$-smoothness of $\bU$
and by the $\psi$-domination of the derivatives. Let  $(X,B) \in \Gn$. 
W.l.o.g. we may assume that the right hand side of \eqref{taylor} is finite. 
Introducing 
\[
\begin{split} 
\eta_{x,x'} \, &:= \, x-x', \quad 
\vartheta_{x,x'} \, :=  \, \tnxb(x') - \tnxb(x) \quad \text{ for } 
 x,x' \in \en(X)\\
&\text{ and } E_{n,\Ko} (X) \,  := \, \{ xx' \in \en(X):  x-x' \notin \Ko  \} \quad 
 \text{ for } X \in \X 
\end{split} 
\]
we have 
\begin{displaymath}
\begin{split}
&H^{\bU}_{\La_{n}}(\iTnB X)\,  
 + \,  H^{\bU}_{\La_{n}}(\Tnb X) \,  
 - \, 2H^{\bU}_{\La_{n}}(X)\\ 
& \,= \sum_{xx' \in E_{n,\Ko}(X)} 
     [ \bU(\eta_{x,x'}  + \vartheta_{x,x'}\einh) 
 + \bU(\eta_{x,x'} - \vartheta_{x,x'}\einh) - 2\bU(\eta_{x,x'})]\\  
& \, \le \sum_{xx' \in E_{n,\Ko}(X)} \psi(x-x') 
 \, (\tnxb(x) - \tnxb(x'))^2  \quad 
=: \quad f_n(X,B).
\end{split}
\end{displaymath}
In the first step we have used that  for $x-x' \in \Ko$ we have $\vartheta_{x,x'} = 0$. 
In the second step we are allowed to apply the above Taylor estimate 
as for  $x-x' \notin \Ko$ we have  $x-x' + s \einh \notin \Ko$ 
for all $s \in [-\vartheta_{x,x'} ,\vartheta_{x,x'} ]$
by \eqref{greps}.  
The arithmetic-quadratic mean inequality gives 
\[
\begin{split}
\frac 1 3 \Big(&(\tnxb(x) - \tn(|x|)) + (\tn(|x|) -\tn(|x'|)) 
 + (\tn(|x'|)- \tnxb(x')) \Big)^2 \\
&\le \, (\tnxb(x) - \tn(|x|))^2 + (\tn(|x|) -\tn(|x'|))^2 
 + (\tn(|x'|)- \tnxb(x'))^2,
\end{split}
\]
and thus 
\begin{displaymath}
\begin{split}
f_{n}(X,B) \,
&\le \, 6 \sideset{}{^{\neq}}\sum_{x,x' \in X} 
 \psi(x-x') \, ( \tn(|x|) - \tnxb (x) )^2 \\ 
&\qquad + 3 \sideset{}{^{\neq}}\sum_{x,x' \in X}
 1_{\{|x| \le |x'|\}}\, \psi (x-x') \, (\tn(|x|)-\tn(|x'|))^2. 
\end{split}
\end{displaymath}
In the first sum on the right hand side we estimate
\[ 
\begin{split}
(\tn&(|x|) - \tnxb (x) )^2 \, 
\le \, \big(\tn(|x|) - \tn(|\anxb(x)|)\big)^2\\
&\le \, \sum_{x'' \in X} 1_{\{|x| \le |x''|\}} 
  1_{\{x \stackrel{X,B_+}{\longleftrightarrow} x''\}} (\tn(|x|) - \tn(|x''|))^2
\end{split}
\]
using \eqref{tauunten}.  By distinguishing the cases 
$x'' \neq x,x'$ and $x'' = x'$ we thus can estimate $f_n(X,B)$ by 
the sum of the two following expressions:
\begin{equation}\begin{split} \label{sig2}
\Sigma_2(n,X) \, &:= \, 9  \sideset{}{^{\neq}} \sum \limits_{x,x'' \in X} 
  \psi(x-x'')  1_{\{|x|\le |x''|\}} |\tn(|x|-\cK) - \tn(|x''|)|^2,\\
\Sigma_3(n,X,B) \, &:= \, 6 \sideset{}{^{\neq}} \sum \limits_{x,x',x'' \in X}
  1_{\{x \stackrel{X,B_+}{\longleftrightarrow} x''\}}
  \psi(x-x')  1_{\{|x|\le |x''|\}}\\
&\hspace{ 4.5 cm} \times |\tn(|x|-\cK) - \tn(|x''|)|^2.
\end{split} \end{equation}
Inserting these sums into the definition of $\Gn$ in  \eqref{good},
we obtain assertion  \eqref{taylor}. 
Assertion \eqref{dichtenklein} can be proved as in Section 
\ref{secngrosshc} using 
\begin{equation} \begin{split} \label{sig4}
\Sigma_4(n,X) \,  &:= \,  2\tau^2 \sideset{}{}\sum_{x \in X} 
  1_{\{x \in \Lan\}} \frac{q(|x| -R)^2}{Q(n-R)^2},\\
\Sigma_5(n,X,B) \, &:=   \, 2 \cf^2 \sideset{}{^{\neq}} \sum_{x,x' \in X} 
 \sideset{}{} \sum \limits_{x'' \in X} 1_{\Kep}(x-x')
   1_{\{x \stackrel{X,B_+}{\longleftrightarrow} x''\}}  1_{\{|x|\le |x''|\}}\\
&\qquad \hspace{3 cm} \times (\tn(|x|-\cK) -\tn(|x''|))^2.
\end{split} \end{equation}
in the definition \eqref{good} of $\Gn$.  

%========================================================================

\subsection{Set of good configurations: Lemma \ref{lebadsmall}}

The functions $\Sigma_i(n,X,B)$ from the definition of the set of 
good configurations $\Gn$ in \eqref{good} have been specified in 
\eqref{sig1}, \eqref{sig2} and \eqref{sig4}. Using the shorthand 
\[
\Taun(x,x'') \, := \, 1_{\{|x|\le |x''|\}} |\tn(|x|-\cK) - \tn(|x''|)|^2
\]
we have 
\[ \begin{split}
& \Sigma_1  =   4 \cf^2 \sideset{}{^{}}\sum \limits_{x,x'' \in X} 
  1_{\{x \stackrel{X,B_+}{\longleftrightarrow} x''\}}
  \Taun(x,x''),\hspace{ 1.8 cm} 
 \Sigma_2  = 9  \sideset{}{^{\neq}} \sum \limits_{x,x'' \in X} 
  \psi(x-x'') \Taun(x,x''), \\ 
& \Sigma_3  =  6 \sideset{}{^{\neq}} \sum \limits_{x,x',x'' \in X}
  1_{\{x \stackrel{X,B_+}{\longleftrightarrow} x''\}}
   \psi(x-x') \Taun(x,x''), \quad 
 \Sigma_4  =  2\tau^2 \sideset{}{}\sum \limits_{x \in X}  1_{\{x \in \Lan\}}
 \frac{q(|x| -R)^2}{Q(n-R)^2},\\ 
& \Sigma_5  =  2 \cf^2  \sideset{}{^{\neq}} \sum \limits_{x,x' \in X} 
 \sideset{}{} \sum \limits_{x'' \in X} 1_{\Kep}(x-x')
   1_{\{x \stackrel{X,B_+}{\longleftrightarrow} x''\}} \Taun(x,x''). 
\end{split} \]
To estimate these sums we set  $\bn := n + \cK$ and $\bR := R + \cK$ 
and use the assertions \eqref{taugegenqhc} and \eqref{ntoinfhc} of Section 
\ref{secsighc}. As a refinement of \eqref{sqmaxhc}, we note that for 
$x_0,\ldots ,x_m \in \R^2$
\[
\Big| |x_m|-|x_0| + \cK \Big| \, 
\le \,  m \bigvee_{i=1}^m  |x_i-x_{i-1}| + \cK  \, 
\le \,  (m+1)( 1 \vee \cK)\Big( 1 \vee \bigvee_{i =1}^m |x_i-x_{i-1}| \Big), 
\]
\begin{equation} \label{sqmax}
\text{ so } \quad (|x_m|-|x_0| + \cK)^2 \, 
    \le \,  (m+1)^2( 1 \vee \cK^2 ) \bigvee_{i =1}^m (1 \vee |x_i-x_{i-1}|^2 ). 
\end{equation}
For the estimation of the expectations of $\Sigma_i$ we combine 
the ideas from Section \ref{secsighc} and from the proof of 
Lemma \ref{lereich}. Using \eqref{taugegenqhc}, \eqref{sqmaxhc}, 
 \eqref{dec}, \eqref{bpsi},\eqref{pig} and \eqref{intg} we obtain
\[
\int \mu(dX) \int \pin(dB|X)\, \Sigma_i(n,X,B)\,\le \, c_i c(n),  
\]
where $c_i$ are finite constants. By \eqref{ntoinfhc} we find that 
\begin{displaymath}
\int \mu \otimes \pin(d(X,B)) \sum_{i=1}^5 \Sigma_i(n,X,B) \, 
\le \, \frac{\de}{4}
\end{displaymath}
for sufficiently large $n$. Now $\mu \otimes \pin (\Gn^c) \le \de$ follows from 
the high probability of  $\Gna$, the Chebyshev inequality and the definition of 
$\Gn$ in \eqref{good}.

\end{sloppypar}

% %===========================================================================
% %
% %===========================================================================

\renewcommand{\thesection}{}

\setlength{\parindent}{0cm}

\end{document}